\documentclass[11pts]{amsart}
\usepackage{amssymb}
\usepackage{amscd}
\usepackage{graphicx}
\usepackage{mathabx}
\usepackage{pdfpages}

\usepackage[
hyperref,
]{degt}

\addtheorem{Notation}{notation}

\def\Dg:{\endgraf{\bf Dg:\enspace}\ignorespaces}
\def\Le:{\endgraf{\bf Le:\enspace}\ignorespaces}
\def\Fl:{\endgraf{\bf Fl:\enspace}\ignorespaces}

\def\sphere{\Bbb S}

\def\ball{\Bbb B}
\def\tsum{\sqcup}
\newcommand{\Zz}{\mathbb{Z}}

\def\bX{\bar X}
\def\tX{\smash{\tilde X}}
\def\dX{\partial X}
\def\bH{\bar H}

\def\bC{\bar C}
\def\bS{\bar S}

\newcommand{\Ss}{\mathbb{S}}

\def\lk{\operatorname{\ell\mathit{k}}}
\def\vlk{\operatorname{\,\overline{\!\lk}}}
\def\sg{\QOPNAME{sg}}
\def\sgn{\QOPNAME{sign}}
\def\Tub{\operatorname{Tub}}
\def\PD{\operatorname{D}}
\def\dm{^*}
\def\ds{^\vee}
\def\units{^{\times\!}}
\def\Cs{\C\units}

\let\uu\upsilon
\def\invsbl#1{\bigl\{#1\bigr\}}

\def\CV{\Cal V}
\def\CA{\Cal A}
\def\CCA{\CA^\circ}
\def\CH{\Cal H}

\def\F{\Bbb F}

\def\ZQ<#1,#2>{\{#1,#2\}}
\def\be{\bold e}

\def\lkg{\operatorname{lkg}}
\let\lkg\kappa

\def\Ann{\QOPNAME{Ann}}

\def\KK{\Cal Z}

\def\CW{{\sl CW}}

\def\ba{\bold{a}}
\def\bb{\bold{b}}
\def\bc{\bold{c}}
\def\bx{\bold{x}}

\def\Log{\operatorname{Log}}

\let\slope\kappa
\def\rGo#1{\vphantom{\Go}\smash{\Go_{#1}}}

\def\vect{}

\let\defect\delta

\def\sg{\operatorname{sg}}

\def\sone{S^1\sminus1}

\def\link#1{\href{http://katlas.math.toronto.edu/wiki/#1}{#1}}

\def\4{\color{magenta}}
\def\3{\color{red}}
\def\7{\color{red}}

\title{Slopes and signatures of links}

\author{Alex Degtyarev}

\address{%
Department of Mathematics\\
Bilkent University\\
06800 Ankara, Turkey}

\email{degt@fen.bilkent.edu.tr}

\author{Vincent Florens}

\address{%
Laboratoire de Math\'{e}matiques et leurs applications, UMR CNRS 5142\\
Universit\'{e} de Pau et des Pays de l'Adour\\
Avenue de l'Universit\'{e}\\
 BP 1155 64013 Pau Cedex, France
}

\email{vincent.florens@univ-pau.fr}

\author{Ana G.\ Lecuona}

\address{%
Aix Marseille Universit\'{e}, CNRS, Centrale Marseille, I2M, UMR 7373, 13453 Marseille, France}

\email{ana.lecuona@univ-amu.fr}

\thanks{%
The first author was partially supported by the T\"{U}B\DOTaccent{I}TAK grant 116F211.
The second author was partially supported by ANR-17-CE40-0023 (LISA).
The third author was partially supported by GEOR MTM2014-55565}



\begin{document}

\begin{abstract}
We define the \emph{slope} of a colored link in an integral homology
sphere, associated to \emph{admissible} characters on the link group.
Away from a certain singular locus, the slope
is a rational function which can be regarded as a multivariate generalization of the
Kojima--Yamasaki $\eta$-function. It is the ratio of two Conway
potentials, \emph{provided that the latter makes sense}; otherwise, it is a
new invariant.
The slope is responsible for an extra correction term in the signature
formula for the splice of two links, in the previously open exceptional case
where both characters are admissible.
Using a similar construction for a special class of tangles, we formulate
generalized skein relations for the signature.
\end{abstract}

\maketitle

\section{Introduction}\label{S.intro}

In our previous paper~\cite{DFL}, we gave a formula for the
signature of the splice of two colored links in terms of the signatures of the summands.
One exceptional case,
that of the characters vanishing on the components of the links identified by
the splice, was left open. The formula does not work in this exceptional case,
and the study of the defect term was the original goal of the present paper.
To do so, we introduce a new link invariant, called \emph{slope}.
This invariant appears to be interesting in its own right,
and we present several properties and examples.

As another development, we change the paradigm and
extend
the splice formula to all, not necessarily rational, characters (see \cite{Vi} for the extension of the signature function).

\subsection{Slopes}\label{intr.slopes}
The \emph{slope} is an isotopy invariant, defined for an oriented link
$K\cup L$, with a distinguished component~$K$, in an integral homology
sphere~$\sphere$.
Optionally, $L$ may be given a $\mu$-coloring
$\pi_0(L)\onto\{1,\ldots,\mu\}$; then, this coloring should be
respected by characters~$\Go$ below.
Denoting by~$T_C$ a small open tubular neighborhood of a link
$C\subset\sphere$, we consider the inclusion
$$
i\:\partial T_K \into \sphere\sminus T_{K \cup L},
$$
with a particular interest in the map induced in the first homology
with the coefficients twisted by a character
$\omega\:\pi_{1}(\Ss\sminus L)\to \C^{\times}$.
If $\omega$ satisfies the
\emph{admissibility} condition $\omega([K])=1$,
then the vector space $H_1(\partial T_K;\C(\Go))$ is
generated by the meridian~$m$ and preferred longitude~$l$ of~$K$.
If the kernel of
 $$
i_*\:H_1(\partial T_K;\C(\Go))\to H_1(\Ss\sminus T_{K\cup L};\C(\Go))
$$
is generated by a single vector $am+bl$, $[a:b]\in\Cp1(\C)$,
we define the \emph{slope} of $K\cup L$ at $\omega$ \via
$$
(K/L)(\Go):= -\frac{a}{b}\in\C\cup\infty
$$
(see \autoref{def.slope}).
Crucial is the fact that, if $\omega$ is unitary and \emph{nonvanishing}
(\ie, $\omega$
 does not take value $1$ on any meridian of $L$), then $(K/L)(\omega)$ is
 well defined and real, possibly infinite (see \autoref{prop.unit}).
This special case is used in our signature formula
(see \autoref{intr.splice} below).

For a link $K\cup L\subset S^3$,
the slope can be computed from the link diagram, using the
 Wirtinger presentation of the group $\pi_1(S^3 \sminus (K \cup L))$ and
 Fox calculus (see \autoref{s.Fox}).
In general, the slope is not determined by the
 combinatorial data, such as linking numbers, depending more
 deeply on the group
(\cf. \autoref{ex.whitehead} of the Whitehead link and a number of
examples in \autoref{s.rational}).

Admissible nonvanishing characters constitute an algebraic variety
$\CCA(K/L)\subset(\CC^\times\sminus1)^\mu$.
Our first main result, \autoref{th.Torres},
states that, if the Alexander polynomial $\Delta_L$ is
  not zero, the function $K/L$ is rational and finite outside
the zero locus of $\Delta_L$. In fact,
denoting by~$\prime$ the derivative with respect to the variable~$t$
corresponding to~$K$,
we have
\[*
(K/L)(\Go)=-\frac{\nabla'_{K\cup L}(1,\sqrt\Go)}{2\nabla_L(\sqrt\Go)}\in\C\cup\infty,
\]
provided that this ratio makes sense, \ie, is not $0/0$.
(Here, $\nabla$ is the Conway potential; we
have to use~$\nabla$ instead of~$\Delta$
to avoid the sign ambiguity, \cf. \autoref{rem.Delta}).
At the common zeroes of the two polynomials, the slope may still be defined,
but its values are less predictable, which makes this invariant interesting.
Even if $\Delta_L\equiv0$, generically the slope is still a rational function
(possibly, but not necessarily, identically~$\infty$) whose denominator is the
first nonvanishing order $\Delta_{L,r}$ (see \autoref{th.rational}).
Our experiments with the link tables~\cite{KAT}
reveal the independence of $K/L$ of the higher Fitting ideals
(\cf. Examples~\ref{ex.slope.2} and~\ref{ex.slope.3});
moreover, the slope distinguishes some links that are not distinguishable by
the higher Alexander polynomials.

As an indirect consequence of \autoref{th.Torres}, it appears that the slope
can be regarded as
  a multivariate generalization of the Kojima--Yamasaki $\eta$-function
  \cite{KY}, defined originally for two-component links with linking number~$0$
  (see \autoref{coro.eta}). Cochran \cite{Co} showed that the
  $\eta$-function provides a generating function of a sequence of
  $\beta$-\emph{invariants} that corresponds to the Sato--Levine invariants
  \cite{Sa} of successive derivatives of the links. In particular, he proved
  that they are integral lifts of certain Milnor-type $\bar{\mu}$-invariants. In a forthcoming paper \cite{DFL2}, we study the concordance invariance of the slope and give more details on its relation to
  the $\eta$-function.
  We also present a number of alternative methods
of
computing the slope,
including those in terms of Seifert surfaces and $C$-complexes.

\subsection{The splice formula}\label{intr.splice}
Our 
 first  motivation for developing the slope function was an attempt to
extend the signature formula for the splice of two colored links to the exceptional
case $\upsilon'=\upsilon''=1$ that was left open in~\cite{DFL}.
Let $L\subset\sphere$ be the splice of two links
$K^*\cup L^*\subset\sphere^*$, where $*=\prime$ or~$\prime\prime$.
Given a unitary character~$\Go^*$ on~$L^*$, denote $\upsilon^*:=\Go^*([K^*])$,
see~\eqref{eq.upsilon}.
Then, as shown in~\cite{DFL}, for a pair $(\Go',\Go'')$ of
\emph{rational} unitary characters one has
\[*
\sigma_L(\omega',\omega'')=
 \sigma_{K'\cup L'}(\upsilon'',\omega')+
 \sigma_{K''\cup L''}(\upsilon',\omega'')
 +\defect_{\lambda'}(\omega')\defect_{\lambda''}(\omega''),
\]
{\em provided that $(\upsilon',\upsilon'')\ne(1,1)$.}
(See \autoref{th.DFL} and \autoref{s.additivity} for the precise statement
and definition of the defect function~$\Gd$, which is related to the signature
of a generalized Hopf link and depends on the linking numbers only.)
In this paper, we consider the signature as a function on all
(non necessarily rational) characters (see also \cite{CNT})
and
 establish
the splice
formula in this full generality.
We mainly
follow
Viro's approach \cite{Vi} for links in $S^3$,
extending it to integral homology spheres
and
filling in a few details.
.

More importantly, we extend the splice formula to the exceptional case where
$\upsilon'=\upsilon''=1$, \ie, both characters are admissible.
Then, we have well defined slopes $\slope^*:=(K^*\!/L^*)(\Go^*)$, and the
formula reads (see \autoref{th.main})
\[*
\sigma_{L}(\Go',\Go'')
 =\sigma_{L'}(\omega')+\sigma_{L''}(\omega'')
 +\delta_{\lambda'}(\Go') \delta_{\lambda''}(\Go'')
 +\Delta\sigma(\slope',\slope''),
\]
where the correction term
$\Delta\Gs(\slope',\slope'')\in\{0,\pm1,\pm2\}$ depends only on the two
slopes, see~\eqref{eq.Delta.sigma}.
(Observe that $\Delta\Gs$ is the only contribution of the knots~$K^*$ to the
formula.) Geometrically, $\Delta\Gs$ is the sum of two Maslov indices in
Wall's non-additivity theorem~\cite{b:Wa}.

We have similar statements for the \emph{nullity}~$\eta_L$ of the splice:
plain additivity if
$(\upsilon',\upsilon'')\ne(1,1)$ (see \autoref{th.DFL})
and an extra correction term
$\Delta\eta(\slope',\slope'')\in\{-1,0,1,2\}$ as the only
contribution of~$K^*$
in the exceptional case (see \autoref{th.main}).

\subsection{Tangles and skein relations}\label{intr.skein}
The
concept of slope turns out quite fruitful in a number of other
applications. As an example, we consider the slope function for
  tangles in integral homology balls with four
labelled loose ends (this time, no component is distinguished).
The main application of this slope is an analogue
 of the skein relations for the signatures.

For a given $\mu$-colored tangle $T$ in $\ball$, with exterior $\ball_T$,
and a character $\omega$ in $(\CC^\times \sminus 1)^\mu$, one can consider the subspace
\[*
\Ker\bigl[\operatorname{inclusion}_*\:H_1(\partial (\ball_{T}) ;\C(\Go))\to H_1(\ball_{T};\C(\Go))\bigr].
\]
These subspaces were already introduced for colored braids and tangles in
\cite{GG,CC1},
assuming that the characters have finite order,
where certain functorial properties with respect to concatenation
were also established (see also \cite{CC2}).
In the present paper, we restrict ourselves to tangles with four ends
and observe that, for \emph{generic} characters $\Go$,
the kernel has dimension one. As in the case of links,
it is then determined by a
single number $\kappa_T(\Go)\in\C\cup\infty$, which we call \emph{slope}.
We study this slope $\kappa_T$ as a function
defined on a dense Zariski
open set of $(\CC^\times\sminus1)^\mu$ ---including all unitary
characters, at which it takes real values.
Using the technics similar to previous sections, we show that
$\kappa_T$ is
the ratio (assuming that it makes sense) of the Conway
potentials of two auxiliary links obtained,
roughly,
by patching the loose ends of the tangle with
an overcrossing and an undercrossing (see \autoref{th.tangle} for the precise
statement).

We define the \emph{sum}
$T'\tsum T''\subset\sphere$
of two tangles $T'\subset\ball'$, $T''\subset\ball''$,
which is a link in an integral homology sphere.
Then, given a pair of unitary characters~$\Go'$, $\Go''$
properly matching on the boundary, we
can define a character $\Go'\tsum\Go''$ on $T'\tsum T''$ and, hence, speak
about the signature $\Gs_{T'\tsum T''}(\Go'\tsum\Go'')$. By Wall's non-additivity theorem~\cite{b:Wa},
 the signatures of the
three pairwise sums of three
tangles $T^i\subset\ball^i$, $i\in\Z/3$ are related
as follows (see \autoref{th.skein})
  \[*
\sum_{i\in\Z/3}\Gs_{T^{i+1}\tsum T^i}(\Go^{i+1}\tsum\Go^i)
 =\sgn(\kappa^0,\kappa^1,\kappa^2).
\]
 Here, $\sgn(\kappa^0,\kappa^1,\kappa^2)\in\{0,\pm1\}$ is a certain
skew-symmetric function of the slopes of the tangles involved,
see \autoref{cor.as.slope}.
For rational characters,
this formula could alternatively be obtained from \cite[Theorem 1.1]{CC1}
and
elementary properties of the Maslov index.

Applying the results to certain elementary basic tangles (see
\autoref{ex.tau+-}), we obtain a multivariate generalization
(and a four-dimensional proof) of the skein
relations found in~\cite{CF} (see \autoref{cor.skein}).
Note that our version is also a refinement: it disambiguates the right hand
side of the relations when it does not make sense (\ie, becomes $0/0$).
In this case, slopes should be computed by other means.
Furthermore, using the concept of slope,
we bridge the gap between the skein relations found
in~\cite{CF}, in terms of the Conway potentials, and those that could be derived from \cite{CC1}, in
terms of the Maslov index.

\subsection{Contents of the paper}\label{intr.contents}
In Section \ref{S.prelim}, we discuss several unrelated well-known subjects
used throughout the paper, our principal goal being fixing the terminology
and notation. Thus, we introduce (co-)homology with twisted coefficients,
characteristic varieties, sign-determined torsion of a \CW-complex, and
twisted Poincar\'{e} duality. We also discuss Wall's non-additivity
theorem, which is our principal tool in establishing the signature formulas,
and recall the terminology
concerning
 colored links and characters on the link group preserving the coloring.

Section \ref{S.slopes} contains the construction of the slope. We study its
 basic properties and, at the end, establish the rationality of the slope and
 its relation to the Conway potential.

Sections \ref{s.signature} and \ref{S.splice} deal with multivariate
signature: we discuss the extension of the signature function to all, not
necessarily rational, characters and prove its invariance
(Section~\ref{s.signature}), and then we state and prove the splice formula,
both for signature and nullity
(Section~\ref{S.splice}).

Finally, in Section \ref{S.skein}, we construct the slope of a
colored tangle, discuss its rationality and relation to the Conway
potential, and state and prove the generalized skein relations.

\subsection{Acknowledgments}
We would like to thank Ken Baker, Anthony Conway, David Cimasoni,
Stepan Orevkov and Enrico Toffoli for their interest in this project
and the very stimulating discussions.

During the work on this paper, the first author was partially supported by
the Associate Scheme at the Abdus Salam International Centre for Theoretical
Physics (ICTP), Trieste; the second author was invited by the
Thematic Trimester on Invariants in low-dimensional geometry and topology of
the CIMI Excellence Laboratory, which was held in Toulouse, France,
(March--July 2017); the third author acknowledges the invitation of the Isaac
Newton Institue to the program ``Homology theories in low dimensional
topology''.

\section{Preliminaries}\label{S.prelim}

 Throughout the paper, all group actions are \emph{right} and by default,
 all modules are right.
 Matrices act on \emph{row vectors} by \emph{right multiplication}.
 We use the notation $R\units$ to refer to the multiplicative group of a
 commutative ring~$R$.

In this section, $X$ will mainly be a finite \CW-complex with $\pi:=\pi_1(X)$. In \autoref{s.Poincare} and \autoref{s.additivity}, we will consider more specifically
 compact smooth manifolds
with boundary, possibly empty;
then, the \CW-decomposition is given by any smooth triangulation.

\subsection{Twisted (co-)homology}\label{s.twisted}

 Define the \emph{chain complex} of right
$\Z\pi$-modules $C_*(X;\Z\pi)$ as
the complex of cellular chains of
the universal covering $\tilde X\to X$,
 freely generated by (arbitrary) lifts of the cells of $X$. For the lifts, we retain the same
orientation as
for the original cells.
Recall that the group ring $\Z\pi$
is equipped with a canonical involutive anti-automorphism
\[*
*\:\Gl=\sum n_ig_i\mapsto\Gl^*:=\sum n_i g_i\1.
\]
(We use $^*$ instead
of~$\bar{\ }$, reserving the latter for the complex conjugation).
Hence, any right $\Z\pi$-module $M$ gives rise to a left
module~$M\dm$, with the $\Z\pi$-action given by $\Gl m:=m\Gl^*$.

For a right $\Z\pi$-module~$M$ and a left $\Z\pi$-module~$N$, we
consider the complexes
\[*
C^*(X;M):=\Hom_{\Z\pi}(C_*(X;\Z\pi),M),\qquad
C_*(X;N):=C_*(X;\Z\pi)\otimes_{\Z\pi}N
\]
and their homology
\[*
H^*(X;M):=H^*(C^*(X;M)),\qquad H_*(X;N):=H_*(C_*(X;N)).
\]
A similar construction applies to a \CW-pair $(X,A)$, starting with
$C_*(X,A;\Z\pi)$.

An important special case is that of a group ring $M=N=\Z G$,
with the $\Z\pi$-module structure given by an
epimorphism $\Gf\:\pi\onto G$. Then, clearly, $C_*$ and $C^*$ are
merely the
\hbox{(co-)}\penalty0chain complexes of the $G$-covering of~$X$
defined by~$\Gf$.
If $G$ is a free abelian group with a basis $t_1,\ldots,t_\mu$,
then $\Z G=\Z[t_1^{\pm1},\ldots,t_\mu^{\pm1}]$ is the ring of Laurent
polynomials;
it is an integral domain, and we can also consider the field of fractions
$\Q(G):=\Q(t_1,\ldots,t_\mu)$.

Another special case is the ring $\C(\Go)$, which is
the field $\C$ regarded as a
$\Z\pi$-module \via\ a multiplicative character $\Go\:\pi\to\Cs$. If $\Go$
factors through a finite abelian group, $\Go\:\pi\onto G\to\Cs$,
the
homology
$H_*(X;\C(\Go))$ can be identified with the $\Go$-isotypical component
(eigenspace) of the induced representation of~$G$ on the homology
$H_*(X_G;\C)$ of the corresponding $G$-covering $X_G\to X$;
in particular, the latter isotypical component does not depend on the choice
of an intermediate group~$G$.
Alternatively, $\C(\Go)$ can be regarded as a local system on~$X$, and then
$H_*(X;\C(\Go))$ is the
ordinary homology of~$X$ with coefficients in this local
system.
Similar statements hold for cohomology and for the (co-)chain complexes.
We
have an obvious $\C$-linear
isomorphism
\[*
C^*(X;\C(\Go))=\Hom_\C(C_*(X;\C(\Go)),\C)
\]
and, hence, a
``universal coefficient formula''
\[
H^*(\,\cdot\,;\C(\Go))=H_*(\,\cdot\,;\C(\Go))\ds,
\label{eq.uc}
\]
were $\ds$ stands for the dual vector space.
We
fix the notation $\Go\dm$ and $\bar\Go$ for
the \emph{dual} and \emph{conjugate} characters, respectively:
\[*
\Go\dm\:g\mapsto\Go(g)\1,\qquad
\bar\Go\:g\mapsto\overline{\Go(g)}.
\]
Note that $\C(\Go)\dm=\C(\Go\dm)$.
We say that~$\Go$ is a \emph{unitary character} if $\bar\Go=\Go\dm$.

The following simple observations, whose utmost importance was probably first
observed by Viro~\cite{Vi}, are used throughout the paper
without further reference.

\lemma\label{lem.circle}
If $\Go\ne1\:H_1(S^1)\to\Cs$, then $H_*(S^1;\C(\Go))=H^*(S^1;\C(\Go))=0$.
\endlemma

\corollary\label{cor.circle.bundle}
If $X\to B$ is a circle bundle and a character $\Go\:\pi\to\Cs$ restricts to
a nontrivial character on the fibers, then $H_*(X;\C(\Go))=H^*(X;\C(\Go))=0$.
\endcorollary

\corollary\label{cor.torus}
If $X=T^2$ is a torus and $\Go\ne1$,
then $H_*(X;\C(\Go))=H^*(X;\C(\Go))=0$.
\endcorollary

The next corollary, although also straightforward, is used several times in
the paper.
Consider
a $3$-manifold~$X$ decomposed into the union $X=Y\cup Z$ of two compact
submanifolds. Given a character $\Go\:\pi_1(X)\to\C\units$, we will say
that the $Z$-part is \emph{$\Go$-invisible} if
\roster*
\item
$Z$ is a plumbed $3$-manifold, possibly with boundary,
which is assumed to be a disjoint union of tori fibered over the boundary
of the core surface, away from the nodes;
\item
$Y\cap Z$ is a union of whole components of the boundary $\partial Z$;
\item
$\Go$
restricts nontrivially to each fiber (``meridian'')
of each circle bundle
constituting~$Z$.
\endroster
To make the $\Go$-invisible part of~$X$ clearly seen, we will
sometimes use
the notation $X=Y\cup\{Z\}$.
Combining Corollaries~\ref{cor.circle.bundle} and~\ref{cor.torus} with the
Mayer--Vietoris exact sequence, we arrive at the following statement.

\corollary\label{cor.plumbing}
Let $X=Y\cup\{Z\}$ be a $3$-manifold and
$\Go\:\pi_1(X)\to\C\units$ a character
such that $Z$ is $\Go$-invisible.
Then the inclusion $Y\into X$ induces isomorphisms
\[*
H_*(Y;\C(\Go))=H_*(X;\C(\Go)),\qquad
H^*(X;\C(\Go))=H^*(Y;\C(\Go)).
\]
\endcorollary

\subsection{Characteristic varieties\pdfstr{}{ \rm(see~\cite{Libgober:char.var})}}\label{s.char.var}

Let $H$ be a free abelian group with basis
$t_1,\ldots,t_\mu$, and let
$\Gf\:\pi\onto H$ be an epimorphism.
Any multiplicative character $\Go\:\pi\to\Cs$ that factors through $H$ is
determined by the sequence
$(\Go_1,\ldots,\Go_\mu)$, where $\Go_i:=\Go(t_i)\in\Cs$.
 This identifies the group of such characters with the complex
torus $(\Cs)^\mu$.

The \emph{characteristic varieties} $\CV_r(X)$ of~$X$ (related to $\Gf$)
are defined \via
\[*
\CV_r(X):=\bigl\{\Go\in(\Cs)^\mu\bigm|\dim H_1(X;\C(\Go))\ge r\bigr\},
 \quad r\ge0.
\]
These are indeed algebraic varieties in $(\Cs)^\mu$,
which
are nested ($\CV_r\supset\CV_{r+1}$) and
depend on the fundamental group $\pi$ of~$X$ (and $\Gf$) only.
In view of~\eqref{eq.conj} below, each $\CV_r$ is real with respect
to the standard real structure $\Go\mapsto\bar\Go$.

Let $C_*:=C_*(X;\Z H)$.
Since the complexes $C_*\otimes_{\Z H}\C(\Go)$ compute the homology
$H_*(X;\C(\Go))$ and $H_0(X;\C(\Go))=0$ unless $\Go=1$, we have
\[*
\CV_r(X)\sminus1=V_\C(E_{r}(C_1/\Im\partial_1))\sminus1,
 \quad r\ge0,
\]
where $E_s(M)\subset\Z H$ is the $s$-th elementary ideal of a
$\Z H$-module~$M$ and $V_\C(I)$, $I\in\Z H$,
stands for the zero locus of the ideal $I\otimes\C\subset\C H$ in
$(\Cs)^\mu=\operatorname{Spec}(\C H)$.
In this identity, $C_*$ can be replaced with any complex of free
$\Z H$-modules computing the homology $H_{\le1}(X;\Z H)$ (\eg, the one
given by the Fox calculus,
\cf. \autoref{s.Fox} below).
In fact, $C_1/\Im\partial_1=H_1(X,x_0;\Z H)$, where $x_0\in X$ is the
basepoint.
Furthermore, according to~\cite{Libgober:char.var}, there also are
isomorphisms
\[*
\CV_r(X)\sminus1=V_\C(E_{r-1}(H_1(X;\Z H)))\sminus1,
 \quad r\ge1,
\]
which are sometimes taken for the definition of $\CV_r(X)$, which extends to
any finitely generated $\Z H$-module.
In particular, it follows that
\[*
\CV_1(X)\sminus1=V_\C(\Ann H_1(X;\Z H))\sminus1.
\]
(Recall that the ideals $E_0(M)$ and $\Ann M$ have essentially the same
radical.)
The irreducible components of $\CV_r(X)$ of codimension $\le1$ constitute the
zero locus of the
\emph{$(r-1)$-st order}
\[*
\Delta_{X,r-1}:=\gcd E_{r-1}(H_1(X;\Z H))
 =\gcd E_{r}(C_1/\Im\partial_1)\in\Z H;
\]
the $0$-th order $\Delta_X:=\Delta_{X,0}$ is called the
\emph{multivariate Alexander polynomial} of~$X$ (and $\Gf$).

\subsection{Sign-determined torsion of \CW-complexes\pdfstr{}{ \rm(see~\cite{Turaev:torsion})}} \label{s.torsion}
Let $\F$ be a field and
\[*
C_*\:\quad
C_m \overset{\partial_{m-1}}{\relbar\joinrel\longto} C_{m-1}
 \overset{\partial_{m-2}}{\relbar\joinrel\longto}
 \cdots \overset{\partial_0}{\relbar\joinrel\longto} C_0
\]
a finite chain complex
of finite-dimensional
vector
 spaces over $\F$.
Suppose that $C_*$ is based, \ie, each~$C_i$ has a
 distinguished basis~$c_i$. Suppose that $H_i(C_*)$ are also based, $i=1,\dots,m$.
  Let
\[*
\beta_i(C_*) :=\sum_{r \le i} (-1)^{i-r}\dim H_r(C_*), \quad
\gamma_i(C_*):=\sum_{r \le i} (-1)^{i-r}\dim C_r, \quad
\ls|C_*|:= \sum_{i=0}^m \beta_i(C_*) \gamma_i(C_*).
\]
Lift the
distinguished
basis for $H_i(C_*)$ to a sequence $h_i\subset C_i$.
Let $b_i\subset C_i$ be a sequence of vectors such that
$\partial_{i-1}(b_i)$ is a basis of $\Im\partial_{i-1}$.
Then, for each~$i$, the
concatenation
$\partial_i (b_{i+1}) h_i b_i$ is a basis of
$C_i$.
Denoting by $T_i$ the transition matrix
from $\partial_i (b_{i+1}) h_i b_i$ to the distinguished basis~$c_i$,
the \emph{torsion} of $C_*$ is defined as
\[*
\tau(C_*) := (-1)^{|C_{*}|} \prod_{i=0}^m \det T_i^{(-1)^{i+1}}\in\F\units,
\]
The torsion of $C_*$ depends on the given bases in $C_i$ and $H_i(C_*)$, but
it
 does not depend on the choice of $h_i$ and $b_i$.

Fix now an
epimorphism $\Gf\:\pi\onto H$ onto a free abelian group $H$ and consider the
complex $C_*(X;\Q(H))$ of $\Q(H)$-vector spaces.
It is based (see \autoref{s.twisted}), and
the \emph{torsion}
$\tau_\Gf(X)\in\Q(H)/{\pm H}$
of~$X$
is defined as the torsion of
$C_*(X;\Q(H))$,
with the extra convention that $\tau_\Gf(X)=0$ if the
complex is not acyclic, \ie, $H_*(X;\Q(H))\ne0$.
Here, the group $\pm H\subset\Q(H)\units$ acts on $\Q(H)$ by the
multiplication, and the ambiguity is
due to the
non-uniqueness in the choice of the bases;
modulo this ambiguity,
$\tau_\Gf(X)$ is invariant under simple homotopy equivalences and subdivisions.

The sign of the torsion can be refined if
$X$ is \emph{homologically oriented}, \ie, equipped with a distinguished
orientation $v$ of the space
$H_*(X;\mathbb{R})=\bigoplus_{i \geq 0}H_i(X;\mathbb{R})$.
The
orientation and order of the cells in \smash{$\tilde{X}$} induce an orientation
and order of those of~$X$, and thus a distinguished bases for
$C_*(X;\mathbb{R})$. Choose bases $h_i$ of
$H_i(X;\mathbb{R})$ so that the basis $h_0 h_1 \cdots h_{\dim X}$ of
$H_*(X;\mathbb{R})$ is positive with respect to $v$, let
$\tau_0(X)=\pm1$ be the sign of the torsion
$\tau(C_*(X;\mathbb{R}) )\in \mathbb{R}\units$, and set
\[*
\tau_\varphi(X,v)= \tau_0(X) \cdot \tau_\varphi(X).
\]
The \emph{sign-determined torsion}
$\tau_\varphi(X,v)\in\Q(G)/H$ depends only on $(X,v)$ and $\varphi$.
In the forthcoming sections,
 $X$ is a $3$-manifold and its torsion is
a topological and simple homotopy invariant.

\subsection{Poincar\'{e} duality\pdfstr{}{ \rm(see~\cite{Wall:surgery})}}\label{s.Poincare}
Throughout this section, $X$ is a smooth compact connected oriented manifold.
 Then,
according to \cite[\S2]{Wall:surgery}, $(X,\dX)$ is a
\emph{simple Poincar\'{e} pair}
of dimension $n:=\dim X$ in the following sense:
there is a simple chain homotopy equivalence
\[*
\PD_X\:C^*(X;\Z\pi)\to C_{n-*}(X,\dX;\Z\pi),
\]
well defined up to chain homotopy; in
particular, for any integer~$r$ and any right $\Z\pi$-module~$M$, there is a
canonical \emph{Poincar\'{e} duality isomorphism}
\[
\PD_X\:H^r(X;M)\overset\cong\longto H_{n-r}(X,\dX;M\dm).
\label{eq.PD.general}
\]
Furthermore, each connected component of the boundary $\dX$ is a simple
Poincar\'{e} complex of dimension $(n-1)$, and the following diagram commutes
\[
\CD
H^r(X;M)@>i^*>>H^r(\dX;M)\\
@V\PD_XV\cong V@V\PD_{\dX}V\cong V\\
H_{n-r}(X,\dX;M\dm)@>\partial>>H_{n-1-r}(\dX;M\dm)\rlap,
\endCD
\label{eq.PD.dX}
\]
where $i\:\dX\into X$ is the inclusion and, for each component $Y\subset\dX$,
we regard $M$ as a $\Z\pi_1(Y)$-module \via\ the inclusion homomorphism
$\pi_1(Y)\to\pi$. (We ignore the technicality related to the choice of
the basepoints as we will mainly work over the commutative rings of the form
$\Z H_1(X)$.)

If $M=\C(\Go)$ for a character $\Go\:\pi\to\Cs$, then,
in view of~\eqref{eq.uc}, the Poincar\'{e} duality in~$X$ and~$\dX$
can be restated in the form of isomorphisms
\[
\gathered
H_{n-r}(X,\dX;\C(\Go\dm))=H_r(X;\C(\Go))\ds,\\
H_{n-1-r}(\dX;\C(\Go\dm))=H_r(\dX;\C(\Go))\ds,
\endgathered
\label{eq.PD.Go}
\]
and \eqref{eq.PD.dX} means that
the
map
$\partial\:H_{n-r}(X,\dX;\C(\Go\dm))\to H_{n-1-r}(\dX;\C(\Go\dm))$
is the adjoint of the inclusion homomorphism
$i_*\:H_r(\dX;\C(\Go))\to H_r(X;\C(\Go))$.

If $n=2r$ is even and the character $\Go$ is unitary,
the first homomorphism in~\eqref{eq.PD.Go}
composed with the
$\C$-anti-linear isomorphism
\[
H_*(\,\cdot\,;\C(\Go))\to H_*(\,\cdot\,;\C(\bar\Go))
\label{eq.conj}
\]
gives rise to
a sesquilinear form
\[
\circ\:\Ga\otimes\Gb\mapsto\<\Ga,\bar\Gb\>\quad\text{on}\quad H_r(X;\C(\Go)),
\label{eq.intersection.index}
\]
called the \emph{intersection index}; it is Hermitian if $r$ is even and
skew-Hermitian if $r$ is odd.
In the former case, the signature $\sign^\Go(X)$
of~\eqref{eq.intersection.index} is called the \emph{twisted signature
of~$X$}. Certainly, in the special case $\Go=1$ we obtain the ordinary
signature $\sign(X)=\sign^1(X)$. The notion of signature extends also
to open manifolds of the form $\bar X\sminus\partial\bar X$, where
$\bar X$ is compact.

By~\eqref{eq.PD.dX}, the kernel $\Ker\partial$ of the intersection index form is
\[*
H_r(X;\C(\Go))^\perp=\Im\bigl[i_*\:H_r(\dX;\C(\Go))\to H_r(X;\C(\Go))\bigr].
\]
In particular, if $\dX=\varnothing$ or, more generally, $H_*(\dX;\C(\Go))=0$,
this form is nondegenerate.

Consider
now the
case of $n=2r+1$ odd and let
\begin{multline}
Z_r(X;\C(\Go)):=\Ker\bigl[i_*\:H_r(\dX;\C(\Go))\to H_r(X;\C(\Go))\bigr]\\
=\Im\bigl[\partial\:H_{r+1}(X,\dX;\C(\Go))\to H_{r}(\dX;\C(\Go))\bigr],
\label{eq.Z}
\end{multline}
where the equality
follows from the exact sequence of
pair $(X,\dX)$.
Since $\dX$ is a closed manifold,
Poincar\'{e} duality induces a perfect pairing
\[
H_r(\dX;\C(\Go))\otimes H_r(\dX;\C(\Go\dm))\to\C
\label{eq.pairing}
\]
which, composed with~\eqref{eq.conj}, coincides with the intersection index
\[
\circ\:\Ga\otimes\Gb\mapsto\<\Ga,\bar\Gb\>\quad\text{on}\quad H_r(\dX;\C(\Go))
\label{eq.Hermitian}
\]
if $\Go$ is unitary.
Combining these observations, we arrive at the following statement.

\lemma\label{lem.PD}
Given a multiplicative character $\Go\:\pi\to\Cs$,
one has
\[*
Z_r(X;\C(\Go))=Z_r(X;\C(\Go\dm))^\perp
\]
with respect to~\eqref{eq.pairing}.
If
$\Go$ is unitary, then
$Z_r(X;\C(\Go))=Z_r(X;\C(\Go))^\perp$
with respect to~\eqref{eq.Hermitian}.
In particular, in this case one has
$\dim Z_r(X;\C(\Go))=\frac12\dim H_r(\dX;\C(\Go))$.
\endlemma

In the next statement, which is an immediate consequence of \autoref{lem.PD},
we change the notation $X\mapsto W$ and $\dX\mapsto\partial W:=X$.

\corollary\label{cor.sign=0}
Assume that a closed oriented $4k$-manifold~$X$ is the boundary
$\partial W$ of a compact
$(4k+1)$-manifold~$W$.
Then
\[*
\sign^\Go(X)=0
\]
for each unitary character $\Go$ on~$X$ which
extends to a unitary character on~$W$.
\endcorollary

\subsection{Signatures and additivity\pdfstr{}{ \rm(see~\cite{b:Wa})}} \label{s.additivity}
Consider a compact connected oriented $4$-manifold $N$ and assume that
$N:=N_1\cup_{X_0}N_2$, where $N_1$, $N_2$ are manifolds with boundaries
\[*
\partial N_1 \cong X_1 \cup_T -X_0,\qquad
\partial N_2\cong X_0\cup_T -X_2
\]
and $X_0$, $X_1$, $X_2$, in turn, are $3$-manifolds with common boundary
\[*
T:=\partial X_0= \partial X_1= \partial X_2.
\]

Consider the $\C$-vector spaces $A_i:=Z_{1}(X_{i};\C)\subset V:=H_{1}(T;\C)$,
$i=0,1,2$, see~\eqref{eq.Z},
and let
$$
W:=\frac{A_0\cap(A_1 + A_2)}{(A_0\cap A_1)+(A_0\cap A_2)}.
$$
By \autoref{lem.PD}, these spaces are Lagrangian with respect to the
intersection index form~$\circ$ of the closed surface~$T$. Hence, $\circ$
induces a nondegenerate Hermitian form~$f$ on~$W$: it is given by
$$
f(a_{0},a_{0}'):=a_{0}\circ a_{1}',\quad
 \text{where $a_{0}'+a_{1}'+a_{2}'=0$  and $a_{i}'\in A_{i}$}.
$$

\theorem[Wall~\cite{b:Wa}]\label{th.wall}
In the above notation, one has
$$
\sign(N)=\sign(N_{1})+\sign(N_{2})-\sign f.
$$
\endtheorem

\remark
As mentioned in~\cite{b:Wa} (and follows easily from the proof), the
conclusion of \autoref{th.wall} holds as well for the twisted signature.
Pick a unitary character $\Go\:H_1(N)\to\Cs$ and denote by the same letter
the restriction of~$\Go$ to the other spaces involved.
Then, the signature formula reads
$$
\sign^{\omega}(N)=\sign^{\omega}(N_{1})+\sign^{\omega}(N_{2})-\sign f,
$$
where the form $f$ is defined as above, using the Lagrangian subspaces
$Z_{1}(X_{i};\C(\omega))$.
\endremark

In the rest of this section, we discuss various forms of Wall's correction
term $\sign f$ in the important special case where $\dim V=2$. Let
$\Gf\:V\otimes V\to\C$ be a nondegenerate skew-Hermitian form; we abbreviate
$\Gf(a,b)=a\circ b$ and $\Gf(a,a)=a^2$. The form $i\Gf$ is
nondegenerate Hermitian; hence,
$\sign(i\Gf)$ takes values $\pm2$ or~$0$. In the former case, $V$ has no
nontrivial isotropic vectors; thus, we assume that $\sign(i\Gf)=0$. In this
case, $V$ has a \emph{standard symplectic basis}, \ie, a basis $m$, $l$
with the property
that
\[*
m^2=l^2=0,\quad m\circ l=-1.
\]
Fixing such a basis, we can parametrize the directions ($1$-subspaces) in~$V$
by assigning to a subspace $\C(am+bl)$ its \emph{slope}
\[
\kappa:=-\frac{a}{b}\in\C\cup\infty=\Cp1(\C).
\label{eq.abstract.slope}
\]
The following statement is immediate.

\lemma\label{lem.isotropic.direction}
A direction in~$V$ is isotropic if and only if its slope is real\rom:
$\kappa\in\R\cup\infty$.
\endlemma

In other words, isotropic directions constitute a circle $\Cp1(\R)$ in the
sphere $\Cp1(\C)=\Bbb{P}(V)$ of all directions. This circle $\R\cup\infty$
has a canonical orientation, \viz. the one that restricts to the order
on~$\R$. (Since the group $\SL(2,\C)$ of isometries of~$V$ is
connected, this orientation does not depend on the choice of a standard
basis.) The Lagrangian subspaces in~$V$ have dimension~$1$; hence, the
correction term $\sign f$ in Wall's formula becomes a function of three
isotropic directions $\C a_i$, $i=0,1,2$, or, equivalently, three slopes
$\kappa_i\in\R\cup\infty$. We will use the notation
\[*
\sign f=\sgn(a_0,a_1,a_2)=\sgn(\kappa_0,\kappa_1,\kappa_2)\in\{0,\pm1\}.
\]

\lemma\label{lem.Wall}
One has
\[*
\sgn(a_0,a_1,a_2)=\sg\bigl[(a_0\circ a_1)(a_1\circ a_2)(a_2\circ a_0)\bigl].
\]
In particular, the function $\sg(a_0,a_1,a_2)$ is skew-symmetric.
\endlemma

\proof
Note that both sides of the identity in the statement remain unchanged if any
of~$a_i$ is replaced with $\Ga_ia_i$, $\Ga_i\in\C\units$. If $A_1=A_2$, both
sides vanish. Otherwise, the pair $(a_1,a_2)$ can be rescaled to a standard
symplectic basis for $V$, \ie, we can assume that $a_1\circ a_2=-1$. Then $W$
is generated by $a_0=\Gb_1a_1+\Gb_2a_2$ and, assuming $\Gb_1$, $\Gb_2$ real,
one has
$f(a_0,a_0)=-\Gb_1\Gb_2$, which is equal to the product
in the right hand side.
\endproof

\corollary\label{cor.as.order}
One has $\sgn(a_0,a_1,a_2)=0$ if and only if at least two of the three
directions coincide. Otherwise, $\sgn(a_0,a_1,a_2)=1$ if and only if
the cyclic order $(a_0,a_1,a_2)$ agrees with the canonical orientation of the
circle of isotropic directions.
\endcorollary

\corollary\label{cor.as.slope}
Given a triple $\kappa_0,\kappa_1,\kappa_2\in\R\cup\infty$, one has
\begin{gather*}
\sgn(\kappa_0,\kappa_1,\kappa_2)=
 \sg\bigl[(\kappa_0-\kappa_1)(\kappa_1-\kappa_2)(\kappa_2-\kappa_0)\bigr]
 \quad\text{if $\kappa_0,\kappa_1,\kappa_2\ne\infty$},\\
\sgn(\infty,\kappa_1,\kappa_2)=\sgn(\kappa_1,\kappa_2,\infty)
 =\sgn(\kappa_2,\infty,\kappa_1)=\sg(\kappa_2-\kappa_1).
\end{gather*}
To make the last formula valid even if one or both of~$\kappa_{1,2}$
is~$\infty$, we extend
the $\sg$ function \via
\[*
\sg x=\begin{cases}
 \hphantom{-}0,&\text{if $x=0$ or $\infty$},\\
 \hphantom{-}1,&\text{if $x>0$},\\
 -1,&\text{if $x<0$}
\end{cases}
\]
and agree to disambiguate $\infty-\infty$ to~$0$ in its argument.
\endcorollary

In a sense, the last formula in \autoref{cor.as.slope} agrees with the
general expression, which gives us
$\sg\bigl[-\infty^2(\kappa_1-\kappa_2)\bigr]$; but then, it disagrees with
our definition of~$\sg$.
For this reason,
we prefer to disambiguate expressions involving~$\infty$ explicitly.



%

\subsection{Colored links}\label{s.colored}

Typically, given an oriented link $L\subset\sphere$
in an integral homology sphere~$\sphere$, we denote by
$T_L:=\Tub L$ a small open
tubular neighborhood of~$L$ and let $X:=\sphere\sminus T_L$.
For a component $C\subset L$, we denote by $\partial_CX$ the intersection of
$\dX$ with the closure of $T_C$.
The group $H_1(\partial_CX)=H_1(\partial T_C)$
is generated by a meridian~$m_C$ and preferred
longitude~$\ell_C$, \viz. the one unlinked with~$C$; we call $\ell_C$ a
\emph{Seifert longitude}.

The meridian~$m_C$ is oriented so that $m_C\circ\ell_C=1$
\emph{with respect to the
orientation of~$\partial T_C$ induced from~$T_C$}.

The group $H_1(X)$ is the free abelian group generated by the classes $m_C$
of the meridians of the components $C\subset L$, and the coloring gives rise to
an epimorphism
\[
\Gf\:\pi_1(X)\longonto H_1(X)\longonto H:=\bigoplus_{i=1}^\mu\Z t_i
\label{eq.Gf}
\]
sending~$m_C$ to~$t_i$ whenever $C\subset L_i$.
Thus, we usually consider the character torus $(\Cs)^\mu$, confining
ourselves to the characters $\Go=(\Go_1,\ldots,\Go_\mu)$ that factor through~$H$.

\definition \label{def.nonvanishing}
A character
$\Go=(\Go_1,\ldots,\Go_\mu)$ is \emph{nonvanishing} if $\Go_i \ne1$ for all $i$.
\enddefinition

The \emph{characteristic varieties}, \emph{orders}, and
\emph{Alexander polynomial} of the colored link~$L$ are defined as those
of~$X$ and $\Gf$, see \autoref{s.char.var}, and we use the notation
$\CV_r(L):=\CV_r(X)$,
$\Delta_{L,r}:=\Delta_{X,r}$, \etc.

Since $H_*(\dX;\C(\Go))=0$ for any nonvanishing character~$\Go$, one has
$H_*(X,\dX;\C(\Go))=H_*(X;\C(\Go))$; hence, by Poincar\'{e}
duality~\eqref{eq.PD.Go},
the \emph{restricted characteristic variety}
$$\CV_r^\circ(L):=\CV_r(L)\cap(\Cs\sminus1)^\mu$$ is invariant under the
automorphism $\Go\mapsto\Go\dm$.
It follows that, modulo units $(\Z H)\units$, the greatest common divisor
$\Delta_r^\circ$ of the defining ideal of the closure of
$\CV_{r+1}^\circ(L)$ is
invariant under the involutive automorphism
$(t_1,\ldots,t_\mu)\mapsto(t_1\1,\ldots,t_\mu\1)$; hence, so is the order
$\Delta_{L,r}$, which differs from $\Delta_r^\circ$ by a number of factors of
the form $(t_i-1)$.

Let, further, $v$ be the orientation of $H_*(X;\R)$ given by the basis
consisting of $1\in H_0(X;\R)$, the meridians $m_C\in H_1(X;\R)$
of the components $C\subset L$ (in some order), and the classes
$[\partial_CX]\in H_2(X;\R)$ of all but the last component (in the same
order). Then, according to~\cite{Turaev:torsion}, up to units $(\Z H)\units$
one has
\[*
\tau_\Gf(X,v)=\begin{cases}
\Delta_L&\text{if $\mu>1$},\\
\Delta_L/(t_1-1)&\text{if $\mu=1$}.
\end{cases}
\]
Hence, as shown above,
\[*
\tau_\Gf(X,v)(t_1\1,\ldots,t_\mu\1)=
 (-1)^{n} t_1^{\nu_1}\ldots t_\mu^{\nu_\mu}\tau_\Gf(X,v)(t_1,\ldots,t_\mu)
\]
for some integers $\nu_1,\ldots,\nu_\mu$,
where $n$ is the number of components of~$L$.
The \emph{Conway potential function}
of the colored link~$L$ is defined as the symmetric
renormalization of~$\tau_\Gf$:
\[
\nabla_L(t_1,\ldots,t_\mu):=
 -t_1^{\nu_1}\ldots t_\mu^{\nu_\mu}\tau_\Gf(X,v)(t^2_1,\ldots,t^2_\mu).
\label{eq.nabla}
\]
The Conway function of~$L$ is that of the maximal coloring; all others are
obtained from the maximal one by the specialization
$H_1(X)\onto H$, \cf.~\eqref{eq.Gf}.

\definition\label{def.nullity}
The \emph{nullity} of a $\mu$-colored link $L\subset\Ss$ at
a nonvanishing character~$\Go$ is
\[
\eta_L(\Go):= \dim \, H_1(X;\Bbb C (\omega)).
\label{eq.nullity}
\]
We extend $\eta$ to all characters \via\ $\eta_L(\Go):=\eta_{L'}(\Go')$,
where $\Go'$ is obtained from~$\Go$ by removing all components $\Go_i=1$, and
$L'$ is obtained from~$L$ by removing the corresponding components~$L_i$.
Occasionally, we consider also the literal extension given by~\eqref{eq.nullity};
then, it is denoted by~$\tilde\eta_L$.
\enddefinition

\remark\label{rem.extension}
In several other definitions below, we use the same strategy as in
\autoref{def.nullity}, \ie, we define a certain quantity $q_L(\Go)$ for
\emph{nonvanishing} characters and extend it to the whole character torus
by patching the components of~$L$ on
which $\Go$ vanishes. The ``literal'' extension, if any, is then denoted
by~$\tilde q_L$. The principal reason for this approach is
\autoref{cor.plumbing}: if there are too many boundary components with
nontrivial homology, we loose control over the situation
(\cf. \autoref{s.omega=1})
and sometimes cannot
even assert that the quantity in question is well defined.
\endremark

A \emph{$(1,\mu)$-colored link} is a $(1+\mu)$-colored link
$
K\cup L
$ in $\mathbb{S}$
in which the \emph{knot} $K$
is the only component (considered distinguished)
given the distinguished color~$0$.
In addition to the space $X=\sphere\sminus L$ and epimorphism~\eqref{eq.Gf},
we will also consider the complement $\bX:=\sphere\sminus(K\cup L)$
and epimorphism
\[
\bar\Gf\:\pi_1(\bX)\longonto H_1(\bX)\longonto \bH:=\Z t\oplus H
\]
sending the meridian $m:=m_K$ to the generator $t:=t_0$.
We use the notation
\[
\vlk(K,L):=(\Gl_1,\ldots,\Gl_\mu)\in\Z^\mu,\quad
\Gl_i:=\lk(K,L_i),\ i=1,\ldots,\mu
\label{eq.vlk}
\]
for the \emph{linking vector} of $K\cup L$.

\section{Slopes}\label{S.slopes}

In this section we define the main character of this article: the slope
of a link in a homology sphere, which is
a function defined on (a part of) the character torus and taking values
in $\C\cup\infty$. We
start
with the definition of this new invariant and its first
properties. We then show how to compute it \via\ the Fox calculus
and, in some special cases, \via\ closed braids. We continue the section by proving that the slope is mainly a rational function, whose poles are determined by the Alexander invariants of the link. Then, we show that in a Zariski dense open set, it is determined by certain Conway potential functions. Finally, we present a list of examples.

\subsection{Definition and first properties}\label{s.defs}

Consider a $(1,\mu)$-colored link $K\cup L\subset\sphere$.

\definition
A
character $\Go\subset(\Cs)^\mu$ of~$\pi_1(\mathbb{S} \smallsetminus L)$ is \emph{admissible}
if $\Go([K])=1$.
 The variety of admissible characters is denoted
$$
\CA(K/L) = \bigl\{ \Go \in (\Cs)^\mu\bigm|\Go^\Gl=1\bigr\},
$$
where \smash{$\Gl:=\vlk(K,L)$} is the linking vector, see~\eqref{eq.vlk}.
\enddefinition

Note that an admissible character $\omega$ restricts to the trivial character
on $H_1(\partial_K\bX)$.
If $\Gl=0$, then $\CA(K/L)=(\Cs)^\mu$;
otherwise,
leting $N:=\gcd(\Gl)$ and
$\nu:=\Gl/N$, the irreducible over~$\Q$ components of $\CA(K/L)$ are
the zero sets of the cyclotomic polynomials $\Phi_d(\Go^\nu)$, $d\mathrel|N$.
Since all varieties are defined over~$\Z$, for each component
$\CA\subset\CA(K/L)$ and each $r\ge0$, the complement $\CA\sminus\CV_r(L)$
is either empty or dense in~$\CA$.

\example \label{ex.linking} If $L$ has one component with $lk(K,L)=\lambda$, then $\CA(K/L) = \{ 1,  \zeta, \zeta^2,\dots, \zeta^{\lambda-1} \}$, where
 $\zeta$ is any primitive root of $1$ of order $\lambda$.
  If $L$ has two components, with $\overline{lk}(K,L)=(1,\lambda)$, then
   $\CA(K/L) = \{ (\omega^\lambda, \omega^{-1}); \, \omega \in  \Cs  \}$.
\endexample

We will mainly consider the variety of nonvanishing admissible
characters
$$\CCA(K/L):=\CA(K/L)\cap(\Cs\sminus1)^\mu.$$
Let $\Go\in\CCA(K/L)$.
Since $\Go$ is
nonvanishing, we have $H_*(\partial_L\bX;\C(\Go))=0$ and, since $\Go$ is also
admissible,
\[*
H_1(\partial\bX;\C(\Go))=H_1(\partial_K\bX;\C(\Go))=H_1(\partial_K\bX;\C).
\]
The latter isomorphism is canonical up to a multiplicative constant,
commutes with \eqref{eq.conj}, and takes~\eqref{eq.pairing}
and~\eqref{eq.Hermitian} to forms that are congruent to, respectively, the bilinear
and sesquilinear extensions of the intersection index form on
$H_1(\partial_K\bX)$.
Hence, we have well-defined subspaces
\[*
\KK(\Go)=\KK_{K\cup L}(\Go):=Z_1(\bX;\C(\Go))\subset H_1(\partial_K\bX;\C),
\]
see \autoref{s.Poincare}, which have the following properties:
\[*
\KK(\bar\Go)=\overline{\KK(\Go)},\qquad
 \KK(\Go\dm)=\KK(\Go)^\perp.
\]
Let $m,l$ be the meridian and
Seifert
longitude of $K$, forming a basis of $H_1(\partial_K\bX;\C(\Go))=\CC^2$.

\definition\label{def.slope}
Let $\Go\in\CCA(K/L)$ and assume that
$\dim\KK(\Go)=1$, \ie, $\KK(\Go)$ is generated by a single vector
$am+bl$ for some $[a:b]\in\Cp1(\C)$. Then,  the \emph{slope} of $K\cup L$ at~$\Go$ is the quotient
\[*
(K/L)(\Go):=-\frac{a}{b}\in\C\cup\infty.
\]
In agreement with \autoref{rem.extension}, we extend $\KK$ and $K/L$ to all
admissible
characters $\Go\in\CA(K/L)$ by patching the components of~$L$ on which $\Go$
vanishes. The literal extension makes no sense.
\enddefinition

\example[generalized Hopf links]\label{ex.Hopf}
Recall
that a \emph{generalized Hopf link} $H_{m,n}$ is obtained from the
ordinary Hopf link $\bar V\cup\bar U$ by replacing $\bar V$ and $\bar U$ by,
respectively, $m$ and $n$ close parallel copies. Assume $m>0$ and take for~$K$
one of the $V$-components. Then, the slope~$\slope$ of $H_{m,n}$ at any
\emph{nonvanishing} character is
\[*
\slope=\begin{cases}
         0 & \mbox{if $n=0$ or $m>1$}, \\
         \infty & \mbox{otherwise}.
       \end{cases}
\]
Indeed, if $n=0$, then $l_K$ bounds a disk, and if $m>1$, then
$l_K$ is homotopic to the longitude of any other $V$-component, which
vanishes in the twisted homology. However, if $m=1$ and $n>0$, then the
meridian~$m_K$ is homotopic to the longitude of any of the $U$-components.
\endexample

The following statements are immediate consequences of Poincar\'{e} duality.

\proposition\label{prop.symmetry}
If the slope at a character $\Go\in\CA(K/L)$ is well defined,
then so are the slopes at
$\bar\Go$, $\Go\dm$, and $\bar\Go\dm$, and one has
\[*
(K/L)(\Go\dm)=(K/L)(\Go),\qquad
(K/L)(\bar\Go)=(K/L)(\bar\Go\dm)=\overline{(K/L)(\Go)}.
\]
Furthermore, the slope does not change if the orientation of~$K$ is reversed.
\endproposition

\proposition\label{prop.unit}
If $\Go\in\CA(K/L)$ is a unitary character,
then the slope $(K/L)(\Go)$ is well defined and,
moreover, is real\rom: $(K/L)(\Go)\in\R\cup\infty$.
\endproposition

\proof
The slope is well defined due to the last statement in \autoref{lem.PD};
it is real due to \autoref{prop.symmetry}.
\endproof

\proposition\label{prop.finite.slope}
The
slope at a character $\Go\in\CCA(K/L)$ is well defined
if and only if the two inclusion homomorphisms
$H_1(K;\C(\zeta))\to H_1(X;\C(\zeta))$, $\zeta=\Go$ or~$\Go\dm$,
are either both trivial or both nontrivial.
The slope is finite, $(K/L)(\Go)\in\C$,
if and only if the two homomorphisms are both trivial.
\endproposition

\proof
By
the Mayer--Vietoris exact sequence, the homomorphism
$H_1(K;\C(\zeta))\to H_1(X;\C(\zeta))$ is trivial if and only if
$\KK(\zeta)$ contains an element of the form $l+\kappa m$, $\kappa\in\C$
(and, in particular, $\KK(\zeta)\ne0$). Thus, the statement follows from the
duality given by \autoref{lem.PD}.
\endproof

\corollary\label{cor.complement}
The slope is
well defined and finite on
$\CA(K/L)\sminus\CV_1(L)$.
\endcorollary

\remark
As follows from \autoref{prop.finite.slope}, for each character~$\Go$,
the existence of the slope
$(K/L)(\Go)$
and its finiteness depend
only on the conjugacy class realized by $K$ in the fundamental group
$\pi_1(\sphere\sminus L)$.
However, the slope itself (when finite) is a more
subtle invariant of $K\cup L$.
\endremark

\proposition\label{prop.ball}
If $K$ is contained in a ball~$B$ disjoint from~$L$,
then $(K/L)(\Go)=0$ for any character $\Go\in\CA(K/L)$.
\endproposition

\proof
Any admissible character restricts to the trivial character on
$\pi_1(B\sminus K)$, and the image of~$l$ vanishes already in
$H_1(B\sminus K;\C(\Go))=H_1(B\sminus K;\C)$.
\endproof

\subsection{Fox calculus}\label{s.Fox}
We illustrate how the slope can be computed by means
of the Fox calculus from a presentation of the fundamental group
$\pi_1(\bar X)$ of the link complement,
together with the classes $m,l\in\pi_1(\bar X)$ of the meridian and
Seifert
longitude of~$K$. In the case of links in~$S^3$, both pieces of data can be
derived from the link diagram.
Indeed, for the group one can choose the Wirtinger
presentation,
where
meridians are the generators.
For~$l$, we trace a curve $C$ parallel to $K$ and such that $\lk(K,C)=0$;
then, starting from the segment corresponding to the chosen meridian of~$K$ and
moving along~$C$ in the positive direction,
we write down the corresponding generator (or its inverse)
each time when undercrossing positively
(respectively, negatively) the diagram of $K\cup L$.
Thus, let
\[*
m,l\in\pi_1(\bar X)=\bigl<x_1,\ldots,x_p\bigm|r_1,\ldots,r_q\bigr>,
\]
$F:=\<x_1,\ldots,x_p\>$, and let $\Lambda:=\Z H$, where $H$ is the
abelianization of~$F$. Since we consider abelian coverings only, we can
specialize Fox derivatives to maps $\partial/\partial x_i\:F\to\Lambda$.
Consider the complex of $\Lambda$-modules \[*
S_*\:\quad S_2\overset{\partial_1}\longto S_1\overset{\partial_0}\longto
 S_0\longto 0,
\]
where
\[*
S_2=\bigoplus_{i=1}^q\Lambda r_i,\quad
S_1=\bigoplus_{i=1}^p\Lambda dx_i,\quad
S_0=\Lambda
\]
and $dx_i$ stands for a formal generator corresponding to~$x_i$.
The ``differential'' of a word $w\in F$ is
\[*
dw:=\sum_{i=1}^p\frac{\partial w}{\partial x_i}dx_i\in S_1;
\]
then, letting
\[*
\partial_1\:r_i\mapsto dr_i,\qquad
\partial_0\:dx_i\mapsto(\text{the image of $x_i$ in $H\subset\Lambda$}),
\]
we obtain a complex computing the homology $H_{\le1}$ of the $H$-covering
of~$\bar X$.

Now, pick an admissible nonvanishing character $\omega \in  \CCA(K/L)$ and
consider the specialization $S_*(\Go):=S_*\otimes_\Lambda\C(\Go)$.
Then, it is straightforward that\[*
\KK(\Go)=\Ker\bigl[ H_1(\partial_KX;\CC(\omega))= \CC m \oplus \CC l
 \overset{\inj_*}\longto
 S_1(\Go) / \Im \partial_1(\Go)\bigr],
\]
where the inclusion homomorphism $\inj_*$ is the specialization of
$m\mapsto dm$, $l\mapsto dl$. (Note that, by the assumption that
$\omega\in\CCA(K/L)$, this homomorphism lands into $\Ker\partial_0(\Go)$.)
Computing the above kernel in the basis $m,l$, we can also compute the slope,
whenever it is defined.

\begin{example}[the Whitehead link]\label{ex.whitehead}
Consider the $(1,1)$-colored Whitehead link $K \cup L$. Since $\lk(K,L)=0$, we have
$\CCA(K/L)=\sone$.
The standard presentation of
$\pi_1(\bar X)$ (derived from the Wirtinger representation) is
\[*
\pi_1(\bar X)=\<m, m_1,l \,|\,[m,l]=1, l=m_1 m\1m_1\1m m_1\1m\1m_1m \>,
\]
where $m$ and $m_1$ are the meridians of $K$ and $L$, respectively,
and $l$ is a Seifert longitude of~$K$.
We can further specialize $\Lambda$ to the group ring
$\ZZ H_1(\bar X)=\Z[t^{\pm1},t_1^{\pm1}]$, sending
the generators~$m$, $m_1$, and~$l$
to $t$, $t_1$, and~$1$, respectively.
Then, denoting by $x$, $y$ the two relations in the presentation above, we have
\[*
dx= (t-1)dl, \quad
dy= dl - t^{-1}(1-t_1)(1-t_1^{-1})dm - (1-t^{-1})(1-t_1^{-1})dm_1.
\]
The specialization at a character $\omega \in  \CCA(K/L)$ means sending
$t_1\mapsto\omega$ and  $t\mapsto1$, so that the image $\Im\partial_1(\Go)$
is generated by $dl-(1-\omega)(1-\omega^{-1})dm$.
We conclude that
$$ (K/L)(\omega)= (1-\omega)(1-\omega^{-1}).$$
\end{example}

The algorithm outlined in this section and using the Wirtinger representation
of the knot group was implemented in \GAP~\cite{GAP4} and used to compute
uni- and multivariate
slopes of all links with up to eleven crossings (see~\cite{KAT}).
We observed quite a few interesting examples, some of which are mentioned in
\autoref{s.rational}.

\subsection{Closed braids} \label{section.burau}
As another example we consider a $(1,1)$-colored link
$K\cup L\subset S^3$, where $L$ is the closure of a braid $\Gb\in\BG{n}$ and
$K$ is its axis. Since $\lk(K,L)=n$,
nonvanishing admissible characters are $n$-th roots of unity $\Go\ne1$.

By the assumption, $S^3\sminus K$ is fibered over the circle~$S^1$, and the
fibers~$D_t$ are open disks that can all be chosen transversal to~$L$. Let
$D^\circ:=D_0\sminus L$, so that $\BG{n}$ acts on the free group
$\pi_1(D^\circ)$. We fix a geometric basis $\Ga_1,\ldots,\Ga_n$ for
$\pi_1(D^\circ)$ and denote by $\Gs_1,\ldots,\Gs_{n-1}$ the corresponding
Artin generators of~$\BG{n}$.
Let $\deg\:\pi_1(D^\circ)\to\Z$ and $\deg\:\BG{n}\to\Z$ be the
homomorphisms given by the exponent sum with respect to the chosen bases.

Let $\Lambda:=\Z[t^{\pm1}]$ be the ring of integral Laurent polynomials.
Recall that the \emph{\rom(reduced\rom) Burau representation} is the
homomorphism $\BG{n}\to\GL(n-1,\Lambda)$ given by the induced $\BG{n}$-action
on the homology $A_n:=H_1(\tilde D^\circ)$ of the infinite cyclic covering
$\tilde D^\circ\to D^\circ$ corresponding to~$\deg$, regarded as a
$\Lambda$-module \via\ the deck translation.
Algebraically, $A_n=\Ker\deg/[\Ker\deg,\Ker\deg]$, and $t$ acts on the class
$[h]$ of an element $h\in\Ker\deg$ \via\ $[h]\mapsto[\Ga_1h\Ga_1\1]$.
As a $\Lambda$-module, $A_n$ is freely generated by the vectors
\[*
\be_i:=[\Ga_{i+1}\Ga_i\1],\quad i=1,\ldots,n-1;
\]
letting $\be_j=0$ for $j\le0$ or $j\ge n$, one has
\[*
\Gs_i\:\ \be_{i-1}\mapsto\be_{i-1}+t\be_i,\quad
\be_i\mapsto-t\be_i,\quad
\be_{i+1}\mapsto\be_i+\be_{i+1},
\label{eq.Burau.n}
\]
(It is
understood
that $\Gs_i\:\be_j\mapsto\be_j$ whenever $\ls|i-j|>1$.)

Since the action of $\BG{n}$ preserves the degree, for any pair
$\Ga\in\pi_1(D^\circ)$ and $\Gb\in\BG{n}$ we have a well defined projection
$\ZQ<\Ga,\Gb>:=(\Ga\Gb)\cdot\Ga\1\in A_n$.
It has the following simple properties:
\[
\begin{alignedat}2
\ZQ<\Ga,\Gb'\Gb''>&=\ZQ<\Ga,\Gb'>\Gb''+\ZQ<\Ga,\Gb''>\quad
 &&\text{for $\Ga\in\pi_1(D^\circ)$ and $\Gb',\Gb''\in\BG{n}$},\\
\ZQ<\Ga,\Gb\1>&=-\ZQ<\Ga,\Gb>\Gb\1\quad
 &&\text{for $\Ga\in\pi_1(D^\circ)$ and $\Gb\in\BG{n}$}.
\end{alignedat}
\label{eq.ZQ}
\]

Pick $\Go=\xi_r:=\exp(2\pi ir/n)\ne1$.
A standard computation shows that, for $\bX=S^{3}\sminus(K\cup L)$ with $L$ the closure of $\beta$, we have
\[
H_1(\bX;\C(\Go))=\bigl(A_n\otimes_\Lambda\C(\Go)\bigr)/(\Gb-1).
\label{eq.braid.H}
\]
The image $l\in H_1(\bX;\C(\Go))$ of
the longitude of~$K$ is the projection of the class of the element
$\Ga_1\ldots\Ga_n$ homotopic to~$\partial D^\circ$; it is easily computed by
the Reidemeister--Schreier method:
\[
l=\sum_{1\le i\le j\le n-1}\Go^j\be_i.
\label{eq.braid.l}
\]
To compute the class of the meridian~$m$,
we identify $H_1(\bX;\C(\Go))$  with the $\Go$-eigenspace of~$t$ in the
homology of the $n$-fold covering $\tX\to\bX$.
Represent~$m$
by a loop $\Gg\in\pi_1(S^3\sminus K)$, which we identify with its lift
to $\pi_1(\tX)$, and pick an element
$\Ga\in\pi_1(D^\circ)$,
$d:=\deg\Ga\ne0\bmod n$.
The homology class
of $\Gg t^d$ is represented by $\Ga\Gg\Ga\1$. On the other hand,
$\Gg\1\Ga\Gg=\Ga\Gb\in \pi_1(\bX)$. Hence,
$(\Go^d-1)m=\ZQ<\Ga,\Gb>$ and, assuming that $\Go^d\ne1$ and
specializing at $t=\Go$,
\[
m=\frac1{\Go^d-1}\ZQ<\Ga,\Gb>(\Go).
\label{eq.braid.m}
\]
In particular, it follows that \eqref{eq.braid.m}
is independent of~$\Ga$.
Formally, $\ZQ<\Ga,\Gb>$ cannot be computed within the framework of the Burau
representation: it is an additional piece of data. In practice, we
usually let $\Ga=\Ga_1$ and compute $\ZQ<\Ga_1,\Gb>$ inductively,
using~\eqref{eq.ZQ} and the obvious identities
\[
\ZQ<\Ga_1,\Gs_1>=t\be_1,\qquad
\ZQ<\Ga_1,\Gs_i>=0\quad\text{for $i>1$}.
\label{eq.braid.ZQ}
\]

According to \autoref{prop.unit}, in the space
$H_1(\bX;\C(\Go))$ given by~\eqref{eq.braid.H} there is exactly one
nontrivial relation $am+bl=0$
between the vectors $m$, $l$ given by~\eqref{eq.braid.m}
and~\eqref{eq.braid.l}, respectively.
Moreover, this relation can be chosen to have real coefficients.
 (We do not know an algebraic
proof of these facts.)
Thus, we have well defined class functions
\[*
 \Gb\mapsto\lkg_r(\Gb):=-\frac{b}{a}
 =\bigl[(K/L)(\xi_r)\bigr]\1\in\R\cup\infty,\quad
r=1,\ldots,n-1.
\]
(Note that we define $\lkg_r$ as the \emph{inverse} of the slope. This
choice is motivated by the somewhat better properties of this invariant
in the realm of closed braids, see, \eg, \autoref{prop.lkg}.)

Numeric experiments reveal rather irregular behaviour of~$\lkg$: it may take
rational or
irrational values (\eg, $\lkg_r(\Gb)=\frac35-\frac15(\xi_r+\xi_r\1)$ for
$\Gb=\Gs_1\Gs_2\Gs_3\Gs_4\Gs_1\Gs_2\in\BG5$),
and it may take value~$\infty$ (\eg, for $\Gb=\Gs_2\1\Gs_1\in\BG3$).
Below are a few further observations concerning these class functions.

\corollary[of~\eqref{eq.braid.m}]
If $\Gb$ has an invariant element $\Ga\in\pi_1(D^\circ)$ of degree~$d>0$
\rom(\eg, if $\Gb\in\BG{d}\times\BG{n-d}\subset\BG{n}$\rom),
then $\lkg_r(\Gb)=0$ whenever $rd\ne0\bmod n$.
\endcorollary

\proposition\label{prop.lkg}
Let $\Gb\in\BG{n}$ and $p\in\Z$. Then, for any $0<r<n$, one has\rom:
\roster
\item\label{B.p}
$\lkg_r(\Gb^p)=p\lkg_r(\Gb)$\rom;
\item\label{B.Delta}
$\lkg_r(\Gb\Delta^2)=\lkg_r(\Gb)-1$\rom;
\item\label{B.p/n}
$\lkg_r((\Gs_1\ldots\Gs_{n-1})^p)=-p/n$.
\endroster
\rom(Here, $\Delta\in\BG{n}$ is the Garside element, so that $\Delta^2$ is the
generator of the center of~$\BG{n}$\rom).
\endproposition

\proof
For~\iref{B.p},
assume that $p>0$; the proof for $p=-1$ and, hence, for $p<0$ is very similar.
From~\eqref{eq.ZQ},
it follows that
$\ZQ<\Ga,\Gb^p>=\ZQ<\Ga,\Gb>\uu$, where
$\uu:=1+\Gb+\ldots+\Gb^{p-1}$.
On the other hand, $\Gb^p-1=(\Gb-1)\uu$, and there remains to
apply the element $\uu\in\Lambda$ to the relation
\[*
am+bl=0\bmod(\Gb-1),
\]
\cf.~\eqref{eq.braid.H}, computing the value $\lkg_r(\Gb)=-b/a$. Here, $l$ and~$m$
are the elements of $A_n\otimes_\Lambda\C(\xi_r)$
given by \eqref{eq.braid.l} and~\eqref{eq.braid.m} (with $\Go=\xi_r$),
respectively; this determines the action of~$\uu$.
One should also observe that
$l$ is $\BG{n}$-invariant,
as so is $\Gr:=[\partial D^\circ]=\Ga_1\ldots\Ga_n$, whence $l\uu=pl$.

For Statement~\iref{B.Delta}, recall that
the action of $\Delta^2$ on $\pi_1(D^\circ)$ is the conjugation by $\Gr\1$
and, since $\deg\Gr=n$,
this action is
identical on $A_n\otimes_\Lambda\C(\Go)$.
Then, since~$l$ is the image of~$\Gr$, we have
$\ZQ<\Ga_1,\Delta^2>=\Gr\1\Ga_1\Gr\Ga_1\1\mapsto(1-\Go)l$ and,
by~\eqref{eq.braid.m} and~\eqref{eq.ZQ}, the new meridian is $m-l$.

For Statement~\eqref{B.p/n}, just recall that
$\Delta^2=(\Gs_1\ldots\Gs_{n-1})^n$ and use~\iref{B.p} and~\iref{B.Delta}.
\endproof

\subsection{Rationality and Conway functions}\label{s.rational}
In this section, we discuss the rationality of the slope as function of~$\Go$
and its relation to the Conway function. We illustrate the results
with a number of examples,
postponing the proofs till the next two sections.

Throughout the section, $K\cup L\subset\sphere$ is a fixed $(1,\mu)$-colored link.

\theorem[see \autoref{proof.rational}]\label{th.rational}
Pick a component $\CA\subset\CA(K/L)$
and let $r$ be the minimal integer such
that $\Delta_{L,r}|_{\CA}\ne0$, \ie, $\CA\sminus\CV_{r+1}(L)$ is dense
in~$\CA$.
Denote by $R$
the coordinate ring of~$\CA$
and fix a normalization of~$\Delta_{L,r}$.
Then, either
\roster
\item\label{rational.finite}
there exists a unique
polynomial $\Delta_{\CA}\in R$ such that
\[*
(K/L)(\Go)=\frac{\Delta_{\CA}(\Go)}{\Delta_{L,r}(\Go)}
\]
holds for each character $\Go\in\CCA\sminus\CV_{r+1}(L)$,
or
\item\label{rational.infinite}
the slope $(K/L)(\Go)=\infty$ is well defined
and infinite
at each character~$\Go$ in a
certain dense Zariski open subset of $\CA$.
\endroster
Case \rom{\iref{rational.infinite}} cannot occur if $r=0$, \ie, if
$\Delta_L|_\CA\ne0$, \cf. \autoref{th.Torres} below.
\endtheorem

\remark\label{rem.limit}
The slope $(K/L)(\Go)$ at a character
$\Go\in\CA\cap\CV_{r+1}(L)$ does not need
to be given by the rational function in
\autoref{th.rational}\iref{rational.finite}, even if the latter
admits an analytic
continuation through~$\Go$.
\endremark

In the next theorem,
as well as in several other statements below, we need to evaluate the Conway
potential at the radical \smash{$\sqrt\Go:=(\sqrt{\Go_1},\ldots,\sqrt{\Go_\mu})$},
which is not quite well defined. We use the convention that one of the values
of each radical is chosen and used consistently \emph{throughout the whole formula};
then, the nature of the formula guarantees that the result is independent of
the initial choice.

\theorem[see \autoref{proof.Torres}]\label{th.Torres}
For a $(1,\mu)$-colored link $K\cup L\subset\sphere$, denote
\smash{$\displaystyle\nabla':=\frac\partial{\partial t}\nabla_{K\cup L}$}.
Then, for a character $\Go\in\CA(K/L)$, one has
\[*
(K/L)(\Go)=-\frac{\nabla'(1,\sqrt\Go)}{2\nabla_L(\sqrt\Go)}\in\C\cup\infty,
\]
provided that the expression in the right hand side makes sense, \ie,
$\nabla'(1,\sqrt\Go)$ and $\nabla_L(\sqrt\Go)$ do not vanish simultaneously.
In particular, the slope is well defined in this case.
\endtheorem

\autoref{th.Torres} is inconclusive if
$\nabla_L(\sqrt\Go)=\nabla'(1,\sqrt\Go)=0$: just as in the freshman calculus,
the indeterminate form $0/0$ should be resolved by other
means. Note also that, even in the case of $(1,1)$-coloring (univariate
polynomials), l'H\^{o}pital's rule does \emph{not} apply!
We illustrate this phenomenon in Examples~\ref{ex.l'Hopital}
and~\ref{ex.l'Hopital.2} below; in a sense, \cf. also \autoref{ex.Hopf}.

\remark
The mysterious polynomial $\frac12\nabla'$
in the statement
can be understood as follows: if $\Go$ is
admissible, then $\nabla_{K\cup L}(1,\sqrt\Go)=0$, \ie,
$\nabla_{K\cup L}(t,\sqrt\Go)=(t-t\1)R(t)$ for a certain Laurent polynomial
$R\in\C[t^{\pm1}]$, and we substitute $t=1$ to the residual
factor~$R$.
\endremark

\remark\label{rem.Delta}
\autoref{th.Torres} can almost be restated in terms of the Alexander
rather than Conway polynomials, thus avoiding the radicals:
the slope is the only monomial multiple of
the ratio
\[*
\pm\frac{\Delta'_{K\cup L}(1,\Go)}{\Delta_L(\Go)}
\]
satisfying \autoref{prop.symmetry}, \ie, such that
$(K/L)(\Go\dm)=(K/L)(\Go)$.
Unfortunately, this simple description misses one vital bit of information
---the sign!.
\endremark

As an indirect consequence of \autoref{th.Torres}, we have that the slope is
a multivariate generalization of the Kojima-Yamaski
$\eta$-function~\cite{KY}.

\corollary \label{coro.eta}
Let $K \cup L$ be a two component $(1,1)$-colored link such that $\lk(L,K)=0$.
Then, for any $\omega \in \sone$ such that $\Delta_{L}(\omega)\neq 0$, the $\eta$-function and the slope coincide at $\omega$.
\endcorollary

\autoref{coro.eta} follows
 from the formula
in \cite{Jin} (first suggested in \cite[Theorem 1]{KY}),
computing
the $\eta$-function in terms of the Alexander polynomials of $K \cup L$ and $L$.
The main theorem in~\cite{Jin} is stated with the
sign ambiguity
(\cf. our \autoref{rem.Delta}),
and the reader is told that the sign can be determined \via\
Bailey's presentation matrix of the first homology group of the universal
abelian cover of $S^{3}\sminus L$.
It is worth noticing
that the slope is defined in the more
general context of links with non-zero linking number. In the
restricted case of two component links, it is defined at each root of unity of order $lk(K,L)$, see \autoref{ex.linking}.
Besides, it may contain certain extra information at the roots
of $\Delta_L$.
(Neither~\cite{KY} nor~\cite{Jin} suggest any clue on the value
of the $\eta$-function at the zeroes of the denominator.)

\begin{example}[the Whitehead link]\label{ex.slope}
Let
$K \cup L$ be the $(1,1)$-colored Whitehead link. We have
$\nabla_{K \cup L}(t,t_1)= (t-t^{-1})(t_1 - t_1^{-1})$ and
$\nabla_L(t_{1})=1/(t_{1}-t_{1}^{-1})$.
Hence, for any
$\Go\in\C\units$,
$$
(K/L)(\omega_{1})=-(\sqrt{\omega_{1}}-\sqrt{\omega_{1}}^{-1})^2=(1-\omega_{1})(1-\omega_{1}^{-1}),
$$
which agrees with \autoref{ex.whitehead}.
This example illustrates also the independence of the ratio
in \autoref{th.Torres} of the choice of~$\sqrt\Go$.
\end{example}

\begin{example}[equal Alexander polynomials]\label{ex.slope.2}
Let $K \cup L_1 \cup L_2$ and $K' \cup L'_1 \cup L'_2$ be the links \link{L11n353}
and \link{L11n384} (see~\cite{KAT}),
respectively.
Both have $11$ crossings and $3$ components, and their Alexander polynomials are
equal:
$$
\Delta_{K \cup L}= \Delta_{K' \cup L'}= (t_2-1)(t-1)^3(t_1-1), \quad
\Delta_{L}= \Delta_{L'}= 0,\quad
\Delta_{L,1}=\Delta_{L',1}=1,
$$
so that \autoref{th.Torres} is inconclusive.
Since $\vlk(K,L)=\vlk(K',L')=(0,0)$, one has
$\CCA(K/L)=\CCA(K'/L')=(\sone)^2$,
and a direct computation using the link diagrams (\cf.
\autoref{s.Fox})
gives us,
for any $\Go:=(\Go_1,\Go_2)\in(\C\units)^2$,
up to the common normalizing factor $-(\Go_1\Go_2)\1$,
\[*
(K/L)(\omega)=(\omega_1 \omega_2^2 + \omega_1^2 - 4 \omega_1 \omega_2 + \omega_2^2+\omega_1),
\quad
(K'/L')(\omega)=(\omega_1 -1)(\omega_1 \omega_2^2-1).
\]
Thus,
the slope can distinguish links with equal Alexander polynomials.
(Here and in the next example, since the first nonvanishing order
$\Delta_{L,1}$ is identically~$1$,
the slope is given by a Laurent polynomial
on the whole torus $(\C\units)^2$, see \autoref{th.rational}.)
\end{example}

\begin{example}[vanishing Alexander polynomial]\label{ex.slope.3}
Let $K \cup L_1 \cup L_2$ be the link \link{L11n396} in~\cite{KAT} with $11$ crossings and
$3$ components.
Both polynomials
$\Delta_{K \cup L}$ and $\Delta_{L}$ vanish identically,
and hence \autoref{th.Torres} fails.
One has $\vlk(K,L)=(0,0)$ and $\CCA(K/L)=(\C\units)^2$.
A direct computation using the link diagrams (\cf.
\autoref{s.Fox}) gives us,
for any $\Go:=(\Go_1,\Go_2)\in(\C\units)^2$,
\[*
(K/L)(\omega)= -(\omega_1 \omega_2- 1)^2/ \omega_1 \omega_2.
\]
Furthermore, the first non-trivial orders are
$\Delta_{K \cup L,1}=(t-1)(t^2-t+1)$ and $\Delta_{L,1}=1$.
These last two examples suggest that the slope
is independent of the higher order Fitting ideals.
\end{example}

\example[l'H\^{o}pital's rule]\label{ex.l'Hopital}
Consider
the family of
two component
algebraically split
links $K\cup L$ described in the left diagram of
\autoref{f.hopital}.
%
\figure
\centering
  \def\svgwidth{\columnwidth}
    \resizebox{0.85\textwidth}{!}{
\begingroup%
  \makeatletter%
  \providecommand\color[2][]{%
    \errmessage{(Inkscape) Color is used for the text in Inkscape, but the package 'color.sty' is not loaded}%
    \renewcommand\color[2][]{}%
  }%
  \providecommand\transparent[1]{%
    \errmessage{(Inkscape) Transparency is used (non-zero) for the text in Inkscape, but the package 'transparent.sty' is not loaded}%
    \renewcommand\transparent[1]{}%
  }%
  \providecommand\rotatebox[2]{#2}%
  \ifx\svgwidth\undefined%
    \setlength{\unitlength}{464.275bp}%
    \ifx\svgscale\undefined%
      \relax%
    \else%
      \setlength{\unitlength}{\unitlength * \real{\svgscale}}%
    \fi%
  \else%
    \setlength{\unitlength}{\svgwidth}%
  \fi%
  \global\let\svgwidth\undefined%
  \global\let\svgscale\undefined%
  \makeatother%
  \begin{picture}(1,0.45851058)%
    \put(0,0){\includegraphics[width=\unitlength]{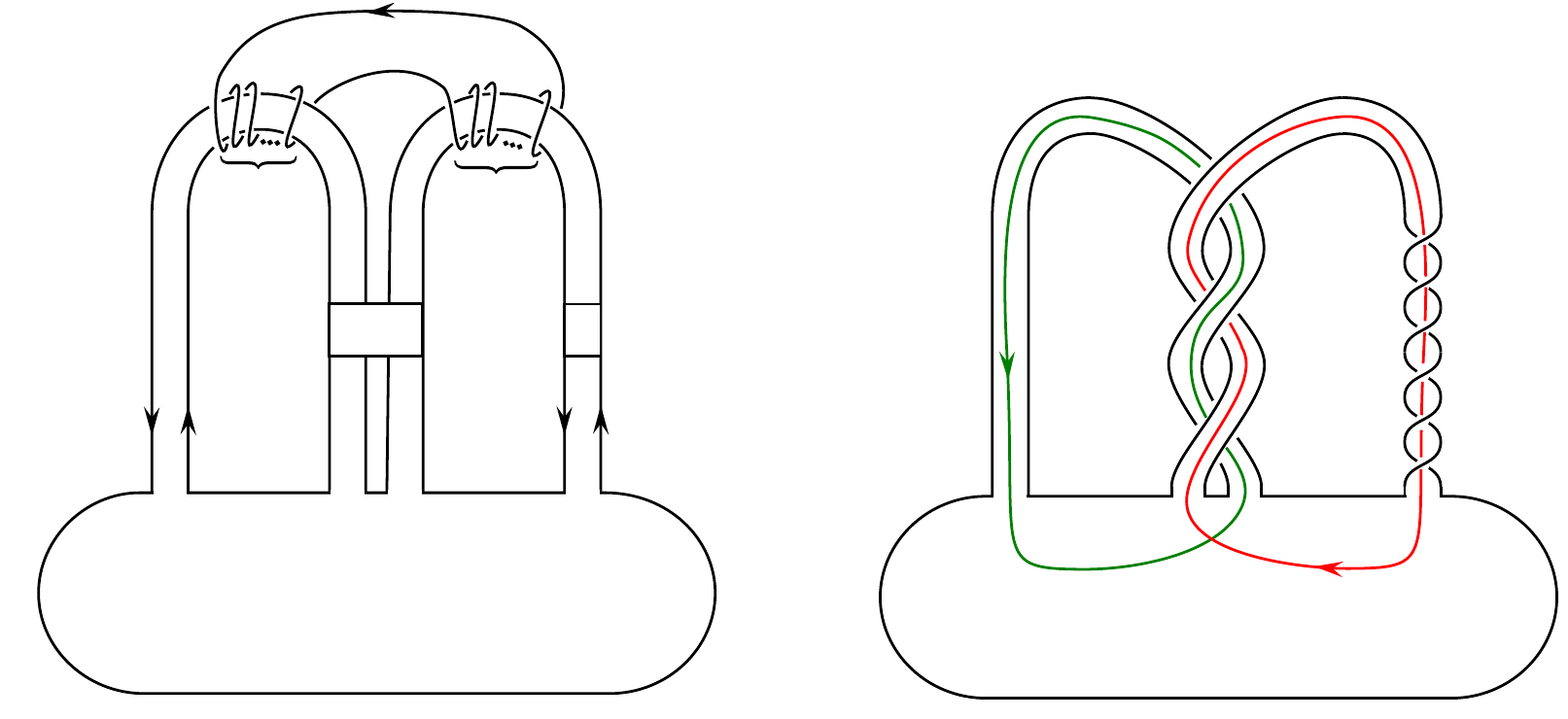}}%
    \put(0.36691763,0.2403089){\makebox(0,0)[lb]{\smash{$c$}}}%
    \put(0.23374003,0.24004888){\makebox(0,0)[lb]{\smash{$a$}}}%
    \put(0.04873374,0.20188393){\color[rgb]{0,0,0}\makebox(0,0)[lb]{\smash{$L$}}}%
    \put(0.12357321,0.43818676){\color[rgb]{0,0,0}\makebox(0,0)[lb]{\smash{$K$}}}%
    \put(0.1575374,0.33222931){\color[rgb]{0,0,0}\makebox(0,0)[lb]{\smash{$x$}}}%
    \put(0.3117244,0.33139893){\color[rgb]{0,0,0}\makebox(0,0)[lb]{\smash{$y$}}}%
    \put(0.61023892,0.02794901){\color[rgb]{0,0,0}\makebox(0,0)[lb]{\smash{$F$}}}%
  \end{picture}%
\endgroup%
}
\caption{ The diagram to the left shows the oriented two component
algebraically split link $K\cup L$ of \autoref{ex.l'Hopital}.
 The one to the right shows
the knot $L$ with parameters $a=3$ and $c=6$.}
\label{f.hopital}
\endfigure %
The knot $L\subset S^{3}$ depends on two
parameters: $a$ stands for an odd number of half crossings between the two
bands, while $c$ stands for an arbitrary number of full crossings between
the two strands. The sign of these two parameters determines whether the
crossings are positive or negative. In the right diagram of
\autoref{f.hopital} the case
$a=3$ and $c=6$ is shown. Setting
$b:=(a-1)/2$,
we obtain that the Seifert matrix for the
Seifert surface~$F$ in \autoref{f.hopital} is given by
\[*
\Theta=\begin{bmatrix}0&b+1\\b&c\end{bmatrix},
\]
and the roots of $\Delta_L$ are  $\Go_\pm:=(1+1/b)^{\pm1}$.
The component $K$ in \autoref{f.hopital}
is the unknot, and the parameters
$x$ and $y$ on the left diagram stand for the linking numbers
between $K$ and the fixed generators for $H_{1}(F)$, depicted in the figure
to the right. In \autoref{f.hopital}, $x\leq 0$ and $y\geq 0$;
changing the direction of the twisting of $K$ around the bands,
one can obtain the other signs. We denote
$\lambda_K:=(x,y)\in\Z^2$.
(In fact, $K$ does not need to be an unknot: the computation
\via\ Seifert
surfaces, which will be explained in~\cite{DFL2}, only makes use of the
linking homomorphism $H_1(F)\to\Z$, $\Ga\mapsto\lk(\Ga,K)$,
\ie, of~$\lambda_K$.)

Computing the slope (see~\cite{DFL2}), we obtain
\[*
(K/L)(\Go)=-\frac{\Go x(2by+y-cx)}{(\Go b-b-1)(\Go b+\Go-b)}
\]
for $\Go\ne1,\Go_\pm$.
Let $\lambda_+:=(2b+1,c)$ and $\lambda_-:=(0,1)$. Further analysis
using~\cite{DFL2} shows that
\roster
\item
if $\lambda_K\ne0,\lambda_\pm$, then
$(K/L)(\Go_\pm)=\infty=\lim_{\Go\to\Go_\pm}(K/L)(\Go)$;
\item
if $\lambda_K=0$, then
$(K/L)(\Go)=0$ for all $\Go\in\Cs\sminus1$;
\item
if $\lambda_K=\lambda_\pm$,
then $(K/L)(\Go)=0$ for all $\Go\in\Cs\sminus\{1,\Go_\pm\}$.
However, in this case we have
$\dim\KK(\Go_\pm)=2$ and $\dim\KK(\Go_\mp)=0$, \ie, the slopes
$(K/L)(\Go_\pm)$ are not defined, even though
$\lim_{\Go\to\Go_\pm}(K/L)(\Go)=0$, (\cf. \autoref{rem.limit}).
\endroster
Now,
we can consider a connected sum~$L'$ of two copies of~$L$, so that the new Seifert
matrix is $\Theta\oplus\Theta$, and choose~$K'$ so that the linking
homomorphism is represented by $\lambda_K:=\lambda_+\oplus\lambda_-$. Then,
$(K'/L')(\Go)=0$ on $\Go\in\Cs\sminus\{1,\Go_\pm\}$
and $(K'/L')(\Go_\pm)=\infty$, whereas
$\lim_{\Go\to\Go_\pm}(K/L)(\Go)=0$.
\endexample

In the previous example, the ``special'' characters $\Go_\pm$ are not
unitary. The next one is more advanced, although slightly less explicit.

\example[l'H\^{o}pital's rule again]\label{ex.l'Hopital.2}
Consider the
two-component link
\link{L10a39} (see~\cite{KAT}). Taking for~$K$ the unknotted component, we have
\[*
\nabla_{K\cup L}=-(t-t\1)(t_1-t_1\1)(t_1^2-1+t_1^{-2})^2,\quad
\nabla_L=(t_1^2-1+t_1^{-2})^2/(t_1-t_1\1),
\]
so that
\[*
(K/L)(\omega)=(\sqrt{\omega}-\sqrt{\omega}^{-1})^2=-(1-\omega)(1-\omega^{-1})
\]
(\cf. \autoref{ex.slope}) unless
$\Delta_L(\Go)=0$, \ie, $\Go=\Go_\pm:=(1\pm i\sqrt3)/2$ is a primitive $6$-th
root of unity. A separate computation using \autoref{s.Fox} shows that
$(K/L)(\omega_\pm)=\infty$ instead of the ``predicted''~$-1$.
Similarly, for \link{L10n36} one has $K/L\equiv0$ except
$(K/L)(\Go_\pm)=\infty$.
\endexample

\example[non-vanishing linking numbers]
Consider the link \link{L10n85} (see~\cite{KAT}) with three (trivial)
components $C_1 \cup C_2 \cup C_3$ ordered and oriented so that
$$
\lk(C_1,C_2)=1,\quad \lk(C_1,C_3)=0,\quad \lk(C_2,C_3)=2.
$$
We have
\[*
\nabla
=(t_2\1-t_2)(t_1^2t_2\1t_3+t_1^{-2}t_2t_3\1
 -t_1^2t_2\1t_3\1-t_1^{-2}t_2t_3
 -2t_2\1t_3-2t_2t_3\1
 +t_2\1t_3\1+t_2t_3)
\]
for the Conway potential of $C_1\cup C_2\cup C_3$.
There are three possible choices for the distinguished component $K$
(where we keep the colors $1,2,3$, making one of them
distinguished and showing the corresponding value $t_i=\bold1$ in bold):
\roster
\item\label{ex.case.1}
If $K=C_1$ and $L=C_2 \cup C_3$, then
$\vlk(K,L)=(1,0)$ and
$\CA(K/L)= \{\bold1\}\times\{ 1 \} \times \C\units$.
 In view of \autoref{rem.extension}, the slope $(K/L)$ is defined by patching
 the component $C_2$, and we have
 $(K/L)(1,\omega)=(C_1/C_3)(\omega)=0$ since $C_1 \cup C_3$ is the trivial link.
\item\label{ex.case.2}
If $K=C_2$ and $L=C_1 \cup C_3$, then $\vlk(K,L)=(1,2)$
and
\[*
\CA(K/L)= \{ (\omega^2,\bold1,\omega^{-1})\,|\,\omega \in \C\units\}.
\]
 We have $\Delta_L=0$ and $\Delta_{L,1}=1$. Moreover,
$\nabla'_2(u,\bold1,u^{-1/2})=2u^{-5/2}(u+1)(u^2-u+1)^2$,
whose roots are
the $3$-rd roots of~$-1$, \ie, $-1$ and $-\xi_3^{\pm1}$.
At $u=1$, the character stops being nonvanishing: $C_1$ needs to be
patched, and the slope of the resulting link \link{L4a1} at $\omega=-1$ is
given by \autoref{prop.lkg}.
The slopes at $u=-\xi_3^{\pm1}$ are computed directly, as in \autoref{s.Fox}.
We obtain
 \[*
\setbox0\hbox{$-2$}
(K/L)(\omega^2,\omega^{-1})=\begin{cases}
         \hbox to\wd0{\hss$0$}, & \mbox{if $\Go=1$}, \\
         \hbox to\wd0{\hss$-2$}, & \mbox{if $\Go^3 = -1$}, \\
         \hbox to\wd0{\hss$\infty$}, & \mbox{otherwise}.
       \end{cases}
\]

\item\label{ex.case.3}
If $K=C_3$ and $L=C_1 \cup C_2$, then $\vlk(K,L)=(0,2)$
and
$\CA(K/L)= \CC\units \times \{ -1,1 \}\times\{\bold1\}$
consists of two components.
For all $\omega \in \CC\units$ we have
$(K/L)(\omega,1)=(C_3/C_1)(\omega)=0$
as in case~\iref{ex.case.1}, whereas
$(K/L)(\omega,-1)=2 (\omega -3 + \omega\1)$
is given by \autoref{th.Torres}, since $L$ is the positive Hopf link and
$\nabla_L=1$ (\cf. also case~\iref{ex.case.2} for the value at $(1,-1)$).
\endroster
 A number of other examples
are given by the
univariate specialization of the slopes of many table links with three or more
components. We have observed all sorts of behaviour of $K/L$ at
unitary roots of $\Delta_L$:
continuity,
infinite value \vs. finite limit, finite value \vs. infinite
limit, \etc.
\endexample

\subsection{Proof of \autoref{th.rational}}\label{proof.rational}
After a sequence of elementary collapses preserving $\partial_K\bX$, we can
assume that $\bX=\sphere\sminus(K\cup L)$ is
a $2$-complex.
Furthermore, the \CW-decomposition of the torus
$\partial_K\bX$ can be chosen standard, with
 a single $0$-cell~$e^0$, two $1$-cells $m$, $l$
 representing the meridian and longitude, respectively,
and a $2$-cell~$k^2$.
Consider the complexes
$$\bS_*:=C_*(\bX;R)\subset S_*:=C_*(X;R).$$
Since they are generated by lifts of the cells to the
corresponding coverings, we retain the same orientation and notation as for
the original cells of $\bX$ and $X$. By the construction of~$\CA$ and~$R$,
both $m$ and $l$ are cycles in~$\bS_*$. Furthermore, $S_1=\bS_1$ and the
image $\partial_1S_2$ differs from $\partial_1\bS_2$
by an extra generator $m=\partial_1e^2$ (for some
extra $2$-cell~$e^2$). Thus, we have the class
\[*
[l]\in H_1(X;R)\subset H_1(X,x_0;R)=\bS_1/(\partial_1\bS_2+Rm),
\]
where $x_0\in X$ is the basepoint. The $R$-module $H_1(X,x_0;R)$
gives rise to a coherent sheaf $\CH$ on~$\CA$, so that
$H_1(X,x_0;R)=\Gamma(\CA;\CH)$,
which restricts to a locally free sheaf (\emph{aka} vector bundle) over
$\CA\sminus\CV_{r+1}(L)$:
the fiber over~$\Go$ is
$H_1(X,x_0;\C(\Go))$, which has constant rank $(r+1)$.
Hence, a section~$s$ of~$\CH$ is in the torsion submodule
$\Tors_RH_1(X,x_0;R)$ if
and only if the support of~$s$ is contained in $\CA\cap\CV_{r+1}(L)$.

We consider separately two cases.

{\em Case \rom{\iref{rational.finite}:} $[l]\in \Tors_RH_1(X,x_0;R)$}.
For any character $\Go\in\CCA\sminus\CV_{r+1}(L)$, there is a polynomial
$p\in R$ such that $p(\Go)\ne0$ and $p[l]=0$ in $H_1(X,x_0;R)$.
This means that
$pl=\tilde pm\bmod\partial_1\bS_2$ for some $\tilde p\in R$ and,
specializing at~$\Go$, we
obtain a
nontrivial
relation $p(\Go)l=\tilde p(\Go)m$
in $H_1(\bX;\C(\Go))$; in particular, $\dim\KK(\Go)\ge1$.
Repeating this argument for $\Go\dm$ (recall that $\CV_{r+1}(L)$ is
symmetric, at least away from the divisors $\Go_i=1$),
we conclude that the slope at~$\Go$ is
well defined and the relation above is the only one, yielding
\[*
(K/L)(\Go)=\frac{\tilde p(\Go)}{p(\Go)}.
\]
This identity holds on the dense set $\{\Go\in\CCA\,|\,p(\Go)\ne0\}$; in
particular, $\tilde p\in R$ is uniquely determined by $p$.
Any two rational functions obtained in this way coincide on a dense set; hence, their
reduced forms are equal, and altogether these functions give rise to a
rational function on $\CCA\sminus\CV_{r+1}(L)$. The denominator of this common
fraction must divide the order $\Delta_{L,r}=\gcd E_{r+1}(H_1(X,x_0;\Z H))$.

{\em Case \rom{\iref{rational.infinite}:} $[l]\notin \Tors_RH_1(X,x_0;R)$}.
In this case, the zero set $Z\subset\CA\sminus\CV_{r+1}(L)$
of $[l]$ regarded as a section of the restricted vector
bundle~$\CH$ is a proper Zariski closed set; by \autoref{prop.finite.slope},
the dense Zariski open subset as in the statement
is $\CCA\sminus(\CV_{r+1}(L)\cup Z\cup Z\dm)$,
where we let $Z\dm:=\{\Go\dm\,|\,\Go\in Z\}$.

The last statement of the theorem follows from \autoref{cor.complement}.
\qed

\subsection{Proof of \autoref{th.Torres}}\label{proof.Torres}
As in the proof of \autoref{th.rational}, after elementary collapses
preserving $\partial_K\bX$, we assume that $\bar X$ is a $2$-complex and that
the \CW-decomposition of
$\partial_K\bX$ has a single $0$-cell~$e^0$, two $1$-cells $m$, $l$
 representing the meridian and longitude, respectively,
and a $2$-cell~$k^2$.
We also assume that $e^0$ is the only $0$-cell of $\bar X$ and that the meridian
$m_C$ of one fixed component $C\subset L$ is represented by a $1$-cell $m_i$.
Denoting by~$r$ the number of $2$-cells, the number of $1$-cells is $r+1$.
(Recall that $\chi(\bX)=0$.)
 We require that
the lift of each $1$-cell~$e^1$ starts at the chosen lift of~$e^0$, and we
order the resulting bases of $\bC_*:=C_*(\bar{X};\ZZ \bar H)$ as follows:
\[
\{e^0\}\subset\bC_0,\quad
\{m,l,\ldots,m_i\}\subset\bC_1,\quad
\{k^2,\ldots\}\subset\bC_2.
\label{eq.bb2}
\]
The same vectors form a basis for $\bS_*:=\bC_*\otimes_{\Z\bH}\Q(\bH)$.

Up to the same simple homotopy equivalence,
the space $X$ is obtained from~$\bX$ by
adjoining an extra $2$-cell $e^2$ bounded by~$m$ and an extra $3$-cell $e^3$
filling the
torus $T_K$.
Besides, the passage from $\Z\bH$ to $\Z H$
is the specialization of the coefficients at
$t=1$.
To respect the order, the generator $e^2$
is to be inserted right after
$k^2$; thus, the distinguished bases for $C_*:=C_*(X;\Z H)$ are
\[
\{e^0\}\subset C_0,\quad
\{m,l,\ldots,m_i\}\subset C_1,\quad
\{k^2,e^2,\ldots\}\subset C_2,\quad
\{e^3\}\subset C_3.
\label{eq.b2}
\]
We orient and lift~$e^2$ and~$e^3$ so that
$\partial_2 e^3=k^2$ and $\partial_1 e^2=m$.

First, assume that $\Go\in\CCA(K/L)$.
Let $b_2:=\{k^2,\ldots\}$ be the chosen basis for $\bS_2$, see~\eqref{eq.bb2}.
With appropriate
orientation and lift of~$k^2$, the matrix of $\partial_1$ has the form
\[
\bmatrix
1-\bar\Gf(l)&t-1&\bold0&0\\
\ba&\bb&M&\bc
\endbmatrix,
\label{eq.d1}
\]
where $M$ is a certain $(r-1) \times(r-2)$-matrix,
$\ba$, $\bb$, $\bc$ are certain
column vectors of dimension~$r-1$, and $\bold0$ is the trivial row vector of
dimension $(r-2)$.
The vector $\partial_0(m_i)=(t_i-1)e^0\ne0$ forms a
basis for~$\bS_0$. It follows that the complex $\bS_*$ is acyclic if and only
if $\partial_1(b_2)m_i$ is a basis for~$\bS_1$, \ie, if the determinant~$D_1$
of the matrix
\[
T_1:=\bmatrix
1-\bar\Gf(l)&t-1&\bold0&0\\
\ba&\bb&M&\bc\\
0&0&\bold0&1
\endbmatrix
\label{eq.T1}
\]
does not vanish. Then, letting $b_0=\varnothing$ and $b_1=\{m_i\}$ in
\autoref{s.torsion}, we obtain
\[*
\tau_{\bar\Gf}(\bX)(t,t_1,\ldots,t_\mu)=
 \frac{-(t-1)\det[\ba\,|\,M]-(\bar\Gf(l)-1)\det[\bb\,|\,M]}{t_i -1}.
\]
Specializing this at $(t,\Go)$ (and thus letting
$\bar\Gf(l)=1$), we arrive at
\[*
\tau_{\bar\Gf}(\bX)(t,\Go)=-\frac{(t-1)\det[\ba\,|\,M](t,\Go)}{\Go_i-1},\quad
\tau'_{\bar\Gf}(\bX)(1,\Go)=-\frac{\det[\ba\,|\,M](1,\Go)}{\Go_i-1}.
\]

A similar computation,
using the bases~\eqref{eq.b2} and matrices
\[*
\partial_1\:\bmatrix
0&0&\bold0&0\\
1&0&\bold0&0\\
\ba&\bb&M&\bc
\endbmatrix_{t=1},\quad
T_1=\bmatrix
1&0&\bold0&0\\
\ba&\bb&M&\bc\\
0&0&\bold0&1
\endbmatrix_{t=1}
\]
instead of~\eqref{eq.d1} and~\eqref{eq.T1}
gives us
\[*
\tau_\Gf(X)(\Go)=\frac{\det[\bb\,|\,M](1,\Go)}{\Go_i-1},
\]
no matter whether $C_*(X;\Q(H))$ is acyclic or not.

Now, in order to compute the slope, we consider the complex $C_*(\bX;\C(\Go))$,
which is merely $C_*(\bX;\Z\bH)$ specialized
at $(1,\Go)$;
we assume this specialization till the end of the computation.
Clearly, a linear combination $am+bl$ represents a class in $\KK(\Go)$ if and
only if it is in the image of~$\partial_1$, which is given by~\eqref{eq.d1},
\ie, essentially, by $[\ba\,|\,\bb\,|\,M\,|\,\bc]$, as the first row
vanishes under the specialization.
For a finite slope~$\kappa$,
we need $-\kappa m+ l\in\Im\partial_1$; clearly,
$\kappa=-\bx\cdot\ba$, where $\bx$ is a solution to the overdetermined
linear system
\[*
\bx\cdot[\bb\,|\,M\,|\,\bc]=[1\,|\,\bold0\,|\,0].
\]
If $\tau_\Gf(X)(\Go)\ne0$, then $\det[\bb\,|\,M]\ne0$ and,
disregarding the last column~$\bc$, we see that
$\KK(\Go)$ contains at most one vector as above, with $\kappa$ given by
Cramer's rule
\[
\kappa
 =-[1\,|\,\bold0]\cdot[\bb\,|\,M]\1\cdot{\ba}
 =-\frac{\det[\ba\,|\,M]}{\det[\bb\,|\,M]}
 =\frac{\tau'_{\bar\Gf}(\bX)(1,\Go)}{\tau_\Gf(X)(\Go)};
\label{eq.Cramer}
\]
in particular, $\dim\KK(\Go)\le1$.
Replacing~$\Go$ with the dual character~$\Go\dm$ and repeating the same
argument, we conclude that also
$\dim\KK(\Go\dm)\le1$; hence, both slopes are well defined and
$\kappa=(K/L)(\Go)$ is given by~\eqref{eq.Cramer}.

If $\tau_\Gf(X)(\Go)=0$ but $\tau'_{\bar\Gf}(\bX)(1,\Go)\ne0$,
\ie, $\det[\bb\,|\,M]=0$ and $\det[\ba\,|\,M]\ne0$, then, arguing
as above and searching for vectors $m-\kappa l\in\KK(\Go)$, we conclude
that $\kappa=0$, \ie, the slope is well defined and
equals~$\infty$.

To complete the proof, we need to take into account the ambiguity
of torsion:
\eqref{eq.Cramer} was obtained by computing both torsions
in compatible distinguished bases.

First, compare the signs (assuming both torsions nonvanishing).
We will use the
notation introduced in \autoref{s.torsion} for~$C_*(X;\R)$ and its barred
counterparts for~$C_*(\bX;\R) \subset C_*(X;\R)$.
Let $n$ be the number of components of~$L$ and $r$, as above, the number of
$2$-cells in~$\bX$.
Clearly,
\[*
b_0=\bar b_0=b_1=\bar b_1=\varnothing,\quad
b_3=\{e^3\}=c_3,\quad
h_0=\bar h_0=\{e^0\}=c_0=\bar c_0,
\]
so that all transition matrices except $T_i$, $\bar T_i$, $i=1,2$, are the
identities.
We can assume that $\bar h_1$ is obtained from~$h_1$ by
prepending $m=\partial_1e^2$. Similarly, we can select
$h_2\subset\bC_2$ and assume that $\bar h_2$
is obtained by prepending $k^2=\partial_2e^3$ to $h_2$, whereas $b_2$ is
obtained by prepending $e^2$ to $\bar b_2$.
Then
\[*
\partial_1(b_2)h_1b_1=(m)\partial_1(\bar b_2)h_1,\qquad
\partial_1(\bar b_2)\bar h_1\bar b_1=\partial_1(\bar b_2)(m)h_1
\]
and, thus,
$\det T_1/\det\bar T_1=(-1)^{\ls|\bar b_2|}=(-1)^{r-n}$,
whereas in dimension~$2$ we have
\[*
\partial_2(b_3)h_2b_2=(k^2)h_2(e^2) \bar b_2,\qquad
\partial_2(\bar b_3)\bar h_2\bar b_2=(k^2)h_2 \bar b_2;
\]
since, on the other hand,
$c_2$ is obtained from~$\bar c_2$ by inserting~$e^2$ as the second
vector, it follows that $\det T_2/\det\bar T_2=(-1)^{\ls|h_2|}=(-1)^{n-1}$.
Since $\ls|C_*(X;\RR)|=1+nr-r$ and $\ls|C_*(\bar X;\RR)|=1+nr$, we conclude that
$\tau_0(X)=-\tau_0(\bX)$, \ie, one should reverse the sign
in~\eqref{eq.Cramer} when switching to the sign-determined torsions.

After passing to the sign-determined torsions, still the quotient in the
right hand side of~\eqref{eq.Cramer} is only well defined up to
multiplicative units~$H$, and there is at most one renormalization of this
quotient
invariant under the involution
$(t_1,\ldots,t_\mu)\mapsto(t_1\1,\ldots,t_\mu\1)$,
see \autoref{prop.symmetry}.
The passage to the Conway functions, through~\eqref{eq.nabla}, is
a way
to obtain such a renormalization
(\cf. also \autoref{rem.Delta}).

Finally,
if $\Go\notin\CCA(K/L)$, we patch the sublinks
$L_i$ corresponding to the vanishing components $\Go_i=1$ and express the
slope in terms of the Conway functions of the two smaller links obtained. The
passage to the Conway functions of the original links is immediate \via\ the
classical Torres relations. (The original article~\cite{To} deals with the
Alexander polynomials, and the translation to the Conway function case
is
found in \cite[Proposition~7]{CiTo}.) Note that, in this last passage, information may be lost, as we may have to multiply both
functions by~$0$.)
\qed

\section{Multivariate signature of colored links} \label{s.signature}

Classically, the
$4$-dimensional approach to the multivariate signature of a colored link uses
branched covers and the $G$-signature theorem  (see, among others, \cite{CF,DFL}).
 Viro \cite{Vi} suggested an alternative
construction, \via\ regular coverings of the complement of the branching
surfaces and cobordisms arguments. This view point
(\cf. also ~\cite{CNT}) allows
one to extend the
signature from rational characters to
the whole character torus $(S^{1}\sminus1)^\mu$.
In this section, we further extend
Viro's construction
to links in integral homology spheres.
At the end, we also deal with the subtleties of
vanishing characters, studying the literal extension of the signature
(\cf. \autoref{rem.extension}) in some special cases.
We
advise
the reader that
most results
of this section
apply to \emph{unitary} characters only.

\subsection{Spanning pairs} \label{s:spanning}
Let $N$ be a compact smooth oriented $4$-manifold with boundary~$\partial N$.
Recall that a compact smooth oriented surface $F\subset N$ is said to be
\emph{properly embedded} if $\partial F=F\cap\partial N\ne\varnothing$ and
$F$ is transversal to~$\partial N$ along $\partial F$.
We define a \emph{properly immersed} surface $F\subset N$ as a finite
union $\bigcup_iF_i$ of connected properly embedded surfaces
$F_i\subset N$ such that all
pairwise intersections of the components~$F_i$ are transversal, at double
points, and away from the boundary $\partial N$.
By a \emph{tubular neighborhood} of $F$ in~$N$ we mean an
\emph{open} regular neighborhood
$T:=T_F\supset F$ which is a union of tubular neighborhoods
$T_i\supset F_i$.

\lemma \label{pair}
Let $F=\bigcup_iF_i\subset N$ be a properly immersed surface.
Fix a tubular neighborhood $T\supset F$ and let
$W:=W_F=N\sminus T$.
Then, the following
three statements are equivalent\rom:
\roster*
\item
$[F_i,\partial F_i]=0\in H_2(N,\partial N;\Q)$ for each index~$i$\rom;
\item
the inclusion homomorphism $H_2(W;\Q)\to H_2(N;\Q)$ is an epimorphism\rom;
\item
the meridians~$m_i$ of all~$F_i$ are linearly independent in $H_1(W)$.
\endroster
Furthermore, the group $H_1(W)$ is generated by the meridians~$m_i$ if and
only if $H_1(N)=0$.
\endlemma

\proof
Let $B$ be the union of all pairwise intersections $T_i\cap T_j$, $i\ne j$:
it is the union of small balls about the points of intersection of the
components of~$F$. Consider $\bar{N}:= N \sminus B$ and let
$\bar F:=F\cap\bar N$ and $\bar T:=T\cap\bar N$. Note that
$\bar F\subset\bar N$
is a properly embedded surface and $\bar T\supset\bar F$ is a tubular
neighborhood in the usual sense. Furthermore, up to homotopy equivalence,
$\bar N$ is obtained from~$N$ by removing a finite set of points, \ie, a
subset of codimension~$4$; therefore,
the inclusion homomorphisms $H_n(\bar N)\to H_n(N)$ are isomorphisms for
$n=0,1,2$, and so are the homomorphisms
$H_2(\bar N,\partial\bar N)\to H_2(N,B\cup\partial N)\leftarrow
H_2(N,\partial N)$.

Consider the exact sequence of the pair
$(\bar N, W)$:
\begin{equation} \label{exact}
 \longto
 H_2(W)\longto
 H_2(\bar N) \stackrel{\rel} \longto
 H_2(\bar N, W) \stackrel\partial\longto H_1(W) \longto H_1(\bar N) \longto
 H_1(\bar N, W) \longto.
\end{equation}
By the excision and Thom isomorphism, we have
\[*
H_n(\bar N,W)=H_n(\bar T,\bar S)=H_{n-2}(\bar F)
\]
(where $\bar S$ is the $S^1$-bundle associated with the disk bundle
$\bar T\to\bar F$).
In particular, the last term in~\eqref{exact} vanishes, whereas the group
$H_2(\bar N, W)$ is generated by the classes $d_i$ of fibers of the disk
bundles $\bar T_i\to\bar F_i$, so that $\partial d_i=m_i$.
This completes the proof of the last statement,
 and the second statement follows immediately.

For the first statement, there remains to observe that the homomorphism
${\rel}\otimes\Q$ is given by
$x\mapsto\sum_i(x\circ[\bar F_i,\partial\bar F_i])d_i$
and, by Poincar\'{e}--Lefschetz duality, ${\rel}\otimes\Q=0$ if and only if all classes
$[\bar F_i,\partial\bar F_i]$ vanish in
$H_2(\bar N,\partial\bar N;\Q)=H_2(N,\partial N;\Q)$.
\endproof

\definition
Let $L$ be a $\mu$-colored link in an integral homology sphere $\Ss$.
A \emph{spanning pair}
for $(\Ss,L)$
is a pair $(N,F)$, where
$N$ is a compact smooth oriented $4$-manifold such that $\partial N=\Ss$ and
$F= F_1 \cup \ldots \cup F_\mu\subset N$ is a properly immersed surface such that
$\partial F_i = F_i \cap \partial N=L_i$ for all $i=1,\dots,\mu$.
We require in addition that $H_1(N)=0$ and
$[F_i,\partial F_i]=0\in H_2(N,\partial N)$ for each index~$i$;
equivalently, we require that the group $H_1(N\sminus F)$
should be freely generated by the
meridians of the components of~$F$.
\enddefinition

The existence of a spanning pair for a colored link is given by
\cite[Proposition 3.4]{DFL}.

\proposition \label{p:boundary}
Fix a $\mu$-colored link $L\subset\Ss$,
consider a spanning pair $(N,F)$, and let $T:=T_F$ be a tubular neighborhood
of~$F$ in~$N$.
Then\rom:
\roster
\item\label{b.character}
each character on~$\Ss\sminus L$ extends to a unique character on $N\sminus T$\rom;
\item\label{b.framing}
for each index~$i$, the Seifert framing of~$L_i$ extends to a framing
of~$F_i$\rom;
\item\label{b.intersection}
for each pair $i\ne j$, the algebraic intersection $F_i\circ F_j$ equals
$\lk(L_i,L_j)$.
\endroster
\endproposition

\proof
Statement~\iref{b.character} is given by \autoref{pair}, since a character, both
on~$L$ and on $N\sminus T$, is uniquely determined by its values on the
meridians. Statements~\iref{b.framing} and~\iref{b.intersection} follow from
the assumption that the classes $[F_i,\partial F_i]$ vanish in
$H_2(N,\partial N)$.
\endproof

\subsection{Invariance of the signature} \label{section.sign}
Fix a $\mu$-colored link $L\subset\Ss$. Given a spanning pair $(N,F)$, we
fix an open tubular neighborhood $T_F$
of~$F$ and let $W_F:=N\sminus  T_F$.

Recall that, according to \autoref{p:boundary}, any character~$\Go$
on~$\Ss\sminus L$
extends to a unique character on $W_F$;
for this reason, we retain the same notation~$\Go$ for
the extension. In this section, we consider unitary characters only.

\definition \label{def.sign}
The \emph{signature} of a $\mu$-colored link $L\subset\Ss$ is the
map
\[*
\sigma_L\:(\sone)^\mu \rightarrow \ZZ,\qquad
 \Go\mapsto \sign^\Go(W_F) - \sign(W_F).
\]
Following \autoref{rem.extension}, we extend the signature function to
arbitrary characters~$\Go$ by patching the components of the link on which
$\Go$ vanishes. Occasionally
(most notably, in the proof of \autoref{th.main}), we need to use the
literal extension~$\tilde\Gs_L$ of the signature, which is not very well
defined; we discuss these subtleties in \autoref{lem.spe} below.

\enddefinition

In view of this definition, in the rest of the paper we
mainly confine ourselves to nonvanishing characters. Furthermore, we usually
use the following alternative definition:
\[
\Gs_L(\Go)=\sign^\Go(W_F)-\sign(N);
\label{eq.W=N}
\]
indeed, by \autoref{pair},
the isometry $H_2(W_F;\Q)\to H_2(N;\Q)$ is surjective and the two forms have
the same signature.

\theorem \label{th.signature}
The signature $\Gs_L$ is independent of the choice of a
spanning pair $(N,F)$.
\endtheorem

\proof
The proof is essentially that of
\cite[Theorem 2.A]{Vi}; we merely fill in a few details.

Given two spanning pairs $(N',F')$, $(N'',F'')$, consider the
closed manifold
$N:=N'\cup_\Ss -N''$ and closed surface $F:=F'\cup -F''\subset N$.
The character~$\Go$ on~$N'$, $N''$ defines a character, also denoted
by~$\Go$, on $W:=N\sminus T_F$.
We can assume that the tubular neighborhoods $T_{F'}\subset N'$ and
$T_{F''}\subset N''$ cut the same tubular neighborhood $T_L\subset\Ss$
of~$L$.
Then, $W=W'\cup_{\Ss\sminus T_L}-W''$ and, by \autoref{th.wall}, we have
\[*
\sign(N)=\sign(N')-\sign(N''),\qquad
\sign^\Go(W)=\sign^\Go(W')-\sign^\Go(W'').
\]
(Indeed, in the former case, the ``corner locus'' $T$ in \autoref{th.wall} is
empty, and in the latter case $T=T_L$ is the union of tori,
each with nontrivial restriction of~$\Go$,
so that
$H_*(T;\C(\Go))=0$ by \autoref{cor.torus}.)
Thus, there remains to prove that $\sign^\Go(W)=\sign(N)$.

By the definition of spanning pair and Mayer--Vietoris exact sequence, we
have $H_1(N)=0$ and $[F_i]=0\in H_2(N)$ for each component~$F_i$ of~$F$.
Since $H_2(X)=\Omega_2(X)$ for any \CW-complex~$X$, each component~$F_i$ is
null-cobordant in~$N$. Pick a cobordism, push it off to
the cylinder $N\times I$, and
smoothen the result to obtain a smooth $3$-manifold $D_i\subset N\times I$
transversal to~$N$ along the boundary $\partial D_i=F_i$.
Do this for each surface~$F_i$ and put the results in general position to
obtain an immersed $3$-manifold $D:=\bigcup_iD_i\subset N\times I$. Let
$T_D\subset N\times I$
be a tubular neighborhood of~$D$, and consider the $5$-manifold
$U:=(N\times I)\sminus T_D$. It is immediate (\cf. \autoref{pair} and the
beginning of this subsection)
that
$H_1(U)$ is generated by the meridians about the components~$D_i$ and,
hence, $\Go$ extends to a
unique character (also denoted by~$\Go$) on~$U$; thus, by
\autoref{cor.sign=0}, we have $\sign^\Go(\partial U)=0$.
On the other hand,
\[*
\partial U=W\cup_{\partial T_F}\partial T_D\sqcup -N
\]
(where $W\subset N\times\{0\}$ and the other copy of~$N$ is $N\times\{1\}$)
and the manifold $\partial T_D$ is obtained by gluing, along whole components
of boundaries, several $4$-manifolds fibered into circles.
Hence, by \autoref{th.wall} and
\autoref{cor.circle.bundle}, we have $\sign^\Go(\partial T_D)=0$ and
\[*
0=\sign^\Go(\partial U)=\sign^\Go(W)-\sign(N),
\]
as stated. (Note that $H_1(N)=0$ and, hence, $\Go=1$ on $N\times\{1\}$.)
\endproof

\subsection{Previous versions of the signature}
Classically,
the signature was defined only for characters \emph{of finite order}, \via\
ramified coverings, and its invariance was proved using the $G$-signature
theorem. We recall briefly the constructions; a more detailed exposition can
be found in~\cite{DFL}.
Let
$(N,F)$ be a spanning pair for $(\mathbb{S},L)$ and
$\Go\in(\Q/\Z)^\mu\subset(S^1)^\mu$
a character of finite order.
Then, $\Go$ defines a normal
covering $N^G \rightarrow N$
with finite abelian group $G\cong\Im\Go$ of deck translations.
Regarding $H_2(N^G;\CC)$ as a $\CC[G]$-module, we consider the eigenspace
\[*
H_2^\omega(N,F):= H_2(N^G;\CC) \otimes_{\CC[G]} \CC
\]
and the restricted
hermitian intersection form; its signature is denoted by $\sign^\omega(N,F)$.
The next lemma asserts that the signature considered in
\autoref{section.sign} extends this definition from the rational points to
the whole character torus
$(\sone)^\mu$.

\begin{lemma} \label{lem.coinc}
For any spanning pair $(N,F)$ and $\omega \in (\sone)^\mu$ of finite order,
one has
$$
\sigma_L(\omega)= \sign^\omega(N,F) - \sign(N).
$$
\end{lemma}

\proof
Let $T_F$ be an open tubular neighborhood of $F$ and $W_F:= N \smallsetminus T_F$.
By \eqref{eq.W=N}, we only need to prove that $\sign^\omega(W_F) = \sign^\omega(N,F)$.
In the notation introduced prior to the statement, we
 have
 an isomorphism of $\CC[G]$-complexes
$$
C_*(W_F; \ZZ \pi_1 (W_F) ) \otimes_{\ZZ \pi_1 (W_F)} \CC[G] \cong  C_*(W_F^G;\CC).
$$
 Hence there is an isomorphism preserving the intersection form
  $$H_2(W_F;\CC(\omega)) \cong H_2^\omega(W_F).$$ The character $\omega$
  induces a branched covering $T_F^G$ of the tubular neighborhood $T_F$,
  branched along $F$. Then $N^G$ is obtained by gluing $T_F^G$ to $W_F^G$,
 along parts of $\partial T_F^G$. The $3$-manifold $\partial T_F$ is a
 plumbing constructed from $F_i \times S^1$, and the fibers
 $\{\cdot\}\times S^1$ are meridians of the components of $L$. Since
 $\omega$ is nonvanishing, by \autoref{cor.circle.bundle},  $T_F^G$ and
 $W_F^G$ are glued along bundles which have trivial homology. By
 Wall's \autoref{th.wall}, we obtain
 $\sign^\omega(W_F) =\sign^\omega(N,F)$.
 \endproof

\subsection{Literal extension of the signature}\label{s.omega=1}
The
following technical lemma, which is used in the proof of
\autoref{th.main},
illustrates the level of difficulties that
one would encounter if the definitions of the nullity and signature were
extended to arbitrary characters literally. Still, we consider a
very special case of one vanishing component only,
and even in this case, the literal extension
$\tilde\Gs_{K \cup L}(1,\omega)$ is not very well defined
unless $\vlk(K,L)\ne0$.

When dealing with spanning pairs of a $(1,\mu)$-colored link
$K\cup L\subset\sphere$, we adopt the notation $(N,D\cup F)$, assuming that
$K=\partial D$ and $L=\partial F$.


\lemma \label{lem.spe}
Let $K \cup L\subset\sphere$ be a $(1,\mu)$-colored link.
Then, for a nonvanishing unitary character
$\Go\in(S^1\sminus1)^\mu$,
one has
\[*
\begin{aligned}
\tilde\sigma_{K \cup L}(1,\omega)&=
\begin{cases}
\sigma_L(\omega)+\sg[(K/L)(\Go)]&
 \text{if $\Go\in\CCA(K/L)$ and $D\cap F=\varnothing$},\\
\sigma_L(\omega) &
 \text{in all other cases},
\end{cases}\\
\tilde\eta_{K \cup L}(1,\omega)&=
\begin{cases}
\eta_L(\omega)+1 &
 \text{if $\Go\in\CCA(K/L)$ and $(K/L)(\omega) \neq \infty$},\\
\eta_L(\omega) &
 \text{in all other cases},
\end{cases}
\end{aligned}
\]
where $\tilde\Gs_{K \cup L}(1,\Go)$ is computed using a
spanning pair $(N,D\cup F)$.
\endlemma


\proof
We can assume that $L\ne\varnothing$, as otherwise both statements become
the tautology $0=0$.

Fix a spanning pair $(N,D\cup F)$.
Consider two transversal tubular neighborhoods
 $T_F\supset F$ and $T_D \cong B^2\times D\supset D$ and introduce
 $W_{D \cup F}:= N \sminus(T_D\cup T_F)$,
 $W_F:= N\sminus T_F$ and
 $W_D:= N\sminus T_D$.
By~\eqref{eq.W=N},
$$
\tilde\sigma_{K \cup L}(1,\omega)= \sign^{1,\omega}(W_{D \cup F}) -
\sign N,
\quad
\sigma_{L}(\omega)= \sign^{\omega}(W_{F}) -
\sign N,
$$
and to prove the first statement of the lemma we will
compare $\sign^{1,\omega}(W_{D \cup F})$ and $\sign^{\omega}(W_{F})$
using Wall's nonadditivity theorem.
The surface $F$ meets~$D$ transversally
in a collection of
$m\ge0$ points,
and
$F \cap T_D$ is a collection of parallel disks
$B_1,\dots,B_m\subset T_D$;
we denote by~$U$ the link $(B_1\cup\ldots\cup B_m)\cap\partial T_D$
in the $3$-manifold $\partial T_D$,
and $T_U\subset\partial T_D$ is its tubular neighborhood
$T_{B_1 \cup \cdots \cup B_m}\cap\partial T_D$.
Let, further, $E^2_m:=D\sminus T_{B_1 \cup \ldots \cup B_m}$, so that
$T_D\sminus T_{B_1 \cup \cdots \cup B_m}\cong B^2\times E^2_m$.
We have
\[*
W_{F} = W_{D \cup F} \cup (B^2\times E^2_m),
\]
glued
along $S^1 \times E^2_m$, where
$S^1 \times \{\cdot\}$ is identified to a meridian of $K$ in $W_{D \cup F}$.

In the rest of the proof, we assume that $\Go\in\CCA(K/L)$ is admissible:
otherwise, the homology $H_1(\partial T_K; \CC(\Go))$ vanishes
(see \autoref{cor.torus}) and the proof
simplifies.

Following the notation of \autoref{th.wall}, we have
$T= \partial T_K \sqcup  \, \partial T_U$
and
\[*
X_0=S^1\times E_m^2,\qquad
X_{1}=(\Ss\sminus T_{K\cup L})\cup\invsbl{\partial T_{F}\cap W_{D}},\qquad
X_{2}=T_{K}\cup\invsbl{\partial T_{F}\cap T_{D}},
\]
where shown in braces $\{\cdot\}$ are $\Go$-invisible parts, see the
definition prior to \autoref{cor.plumbing}.
Since $\omega$ is a nonvanishing character,
we have $H_1(\partial T_U; \CC(\omega))=0$.
Since $\Go$ is also admissible, we also have
$H_1(T;\CC(1,\omega))=H_1(\partial T_K; \CC(1))=\CC^2$.
Thus, using \autoref{cor.plumbing} and ignoring the $\Go$-invisible parts,
we can easily compute the subspaces $A_i\subset H_1(X_i;\C(1,\Go))$ in
\autoref{th.wall}:
\begin{alignat*}2
H_1(X_0;\CC(\omega))&= H_1(S^1\times E_m^2; \CC(1,\omega)),\qquad
 & A_0&=
 \langle m_K \rangle\ \text{or}\ \<l_K\>,\\
H_1(X_1;\CC(1,\omega))&= H_1(\mathbb{S} \sminus T_{K \cup L};\CC(1,\omega)),
 & A_1&= \mathcal{Z}_{K \cup L}(\omega)= \langle a m_K + b l_K \rangle,\\
H_1(X_2; \CC(1,\omega))&= H_1(T_K; \CC(1)),
 & A_2&= \langle m_K \rangle,
\end{alignat*}
where $a,b\in\CC$ and $-a/b=(K/L)(\Go)$.
Here, $A_0=\<m_K\>$ if $m>0$ and $A_0=\<l_K\>$ if $m=0$; this space is
computed exactly as in \autoref{ex.Hopf}. (If $m=0$, we also use the
obvious fact
that the Seifert framing of~$K$ extends to a framing of~$D$,
so that $l_K$ bounds a parallel copy of~$D$.)

If $m>0$, then $A_0=A_2$, the
correction term $\sign f$ in~\autoref{th.wall} vanishes
(just like in the easier case where $\Go$ is not admissible,
where $A_0=A_1=A_2=0$), and we obtain,
in both cases,
\[*
\sign^{\omega}(W_F)=\sign^{1,\omega}(W_{D \cup F})
 +\sign^{\omega}( B^2\times E_m^2).
\]

If $m=0$, the
left hand side of this last expression acquires an extra summand
\[*
\sign f=\sgn(0,\slope,\infty)=\sg\slope,\qquad
\slope:=-a/b=(K/L)(\Go),
\]
see \autoref{cor.as.slope}.
There remains to observe that $E_m^2$ is a surface with
nonempty boundary and, hence,
we have
$H_2(B^2\times E_m^2;\C(\Go))=H_2(E_m^2;\C(\Go))=0$
and $\sign^{\omega}(B^2\times E_m^2)=0$.

The formula for the nullity follows from the Mayer--Vietoris exact
sequence
\begin{multline*}
0 \longto \mathcal{K}
  \longto H_1(\partial T_{K}; \CC(1))
  \longto  \\
   H_1(\Bbb S \smallsetminus (K \cup L);\CC(1,\omega )) \oplus H_1(T_{K};\CC(1))
  \longto H_1(\Bbb S \smallsetminus L; \CC(\Go))
  \longto 0,
\end{multline*}
where
$\mathcal{K}=\mathcal{Z}_{K \cup L}(\omega) \cap Z_1(T_K,\CC(1)){}=A_{1}\cap A_{2}$.
Hence, we have $\dim \mathcal{K}\le1$ and $\dim \mathcal{K}=1$
if and only if $(K/L)(\omega)=\infty$, implying the statement.
(Recall
that we assume the character both
unitary and admissible and, hence, the slope is well defined.)
\endproof

\section{The splice formula} \label{S.splice}

Recall
that the \emph{splice} of two $(1,\mu^*)$-colored links
$K^*\cup L^*\subset\Ss^*$, $*=\prime$ or $\prime\prime$,
is defined as follows.
Denote by
$T^*\subset\Ss^*$ a
small tubular neighborhood of~$K^*$ disjoint from~$L^*$ and let
$m^*,l^*\subset\partial T^*$ be
its meridian and
Seifert longitude, respectively.
Then, the splice of the two links
is the $(\mu'+\mu'')$-colored link $L:=L'\cup L''$ in the integral
homology sphere
$$
\Ss:=(\Ss'\smallsetminus T')\cup_\varphi(\Ss''\smallsetminus T''),
$$
where the gluing homeomorphism $\varphi\colon\partial T'\to\partial T''$
takes~$(m',l')$ to~$(l'',m'')$, respectively.

\subsection{Statement of the splice formula}

A formula for the colored signature of the splice of two links was established in \cite{DFL}, under some restrictions on the characters.
In order to state the general formula, we first introduce some notation.

Let $L\subset\Ss$ be the splice of the $(1,\mu^*)$-colored links
$K^*\cup L^*\subset\Ss^*$, $*=\prime$ or $\prime\prime$.
Consider the linking vectors $\lambda^{*}:=\vlk(K^{*},L^{*})$,
see~\eqref{eq.vlk},
and,
for characters \smash{$\omega^*\in (S^1)^{\mu^*}$}, denote
\[
\upsilon^*:=\Go^*([K^*])
 =\prod_{i=1}^{\mu^{*}}(\Go^{*}_{i})^{\lambda^{*}_{i}}
 =(\omega^*)^{\lambda^*}\in S^1 \subset \CC^\times.
\label{eq.upsilon}
\]
Define the \emph{defect function}
\begin{align*}
\textstyle
\defect_\lambda \colon (S^1)^\mu&\longto\Z \\
 \vect\omega&\longmapsto\textstyle
 \ind\bigl(\sum_{i=1}^\mu \lambda_i \Log\omega_i\bigr)-\sum_{i=1}^\mu \lambda_i \ind(\Log\omega_i),
\end{align*}
where $\vect\lambda\in\Zz^\mu$, the \emph{index} of a real number $x$ is
 defined \via\ $\ind(x):=\lfloor x\rfloor-\lfloor-x\rfloor\in\ZZ$,
 and the \emph{$\Log$-function}
$\Log\colon S^1\to[0,1)$ sends $\exp(2\pi it)$
to $t\in[0,1)$.
Then, the main result in \cite{DFL}
can essentially be stated as follows.

\theorem[\cf. \cite{DFL}]\label{th.DFL}
In the notation introduced at the beginning of the section,
under the assumption that $(\upsilon',\upsilon'')\ne(1,1)$, one has
\[*
\begin{aligned}
\sigma_L(\omega',\omega'')&=
 \sigma_{K'\cup L'}(\upsilon'',\omega')+
 \sigma_{K''\cup L''}(\upsilon',\omega'')
 +\defect_{\lambda'}(\omega')\defect_{\lambda''}(\omega''),\\
\eta_L(\omega',\omega'')&=
 \eta_{K'\cup L'}(\upsilon'',\omega')+
 \eta_{K''\cup L''}(\upsilon',\omega'').
\end{aligned}
\]
\endtheorem

The assumption that $(\upsilon',\upsilon'')\ne(1,1)$ is crucial, as
\cite[Example 2.5]{DFL} shows.

Strictly speaking, the signature formula is proved in~\cite{DFL} only for
rational characters~$\Go^*$, \ie, such that \smash{$\Log\Go^*\in\Q^{\mu^*}$}.
However, once the signature is defined, the extension of the formula to the
whole character torus $(\sone)^\mu$
is immediate, as
Wall's non-additivity formula still works.
(Alternatively, one can follow the proof found in \autoref{section.proof}
below, omitting all slope computations, as the homology groups of all tori
involved vanish.) The nullity formula, not stated explicitly in~\cite{DFL},
follows from \autoref{cor.torus} and the Mayer--Vietoris exact sequence
(again, \cf. the more involved case treated in \autoref{section.proof}).

Our goal is extending \autoref{th.DFL} to the special case
$\upsilon'=\upsilon''=1$.
By definition, $\upsilon^*=1$ if and only if
$\Go^*$ is an admissible character, \ie, there is a well defined slope
$\slope^*:=(K^*\!/L^*)(\Go^*)$.
These slopes give rise to an extra correction term, described
in the following statement.

\theorem[see \autoref{section.proof}]\label{th.main}
Consider two $(1,\mu^*)$-colored links $K^*\cup L^*\subset\Ss^*$
as at the
beginning of the section
and their splice $L:=L'\cup L''\subset\Ss$.
Let \smash{$\Go^*\subset\CA(K^*\!/L^*)\cap(S^1)^{\mu^*}\!$} be two admissible
characters
\rom(so that $\upsilon'=\upsilon''=1$, see~\eqref{eq.upsilon}\rom),
and denote
$\slope^*=(K^*\!/L^*)(\omega^*)$. Then
\[*
\begin{aligned}
\sigma_{L}(\Go',\Go'')&
 =\sigma_{L'}(\omega')+\sigma_{L''}(\omega'')
 +\delta_{\lambda'}(\Go') \delta_{\lambda''}(\Go'')
 +\Delta\sigma(\slope',\slope''),\\
\eta_{L}(\Go',\Go'')&
 =\eta_{L'}(\omega')+\eta_{L''}(\omega'')
 +\Delta\eta(\slope',\slope''),
\end{aligned}
\]
where the correction terms
$\Delta\Gs,\Delta\eta\in\{0,\pm1,\pm2\}$
are as shown in \autoref{fig.main}
\rom(see \autoref{rem.maim}\rom).
\figure
\includegraphics[width=11cm]{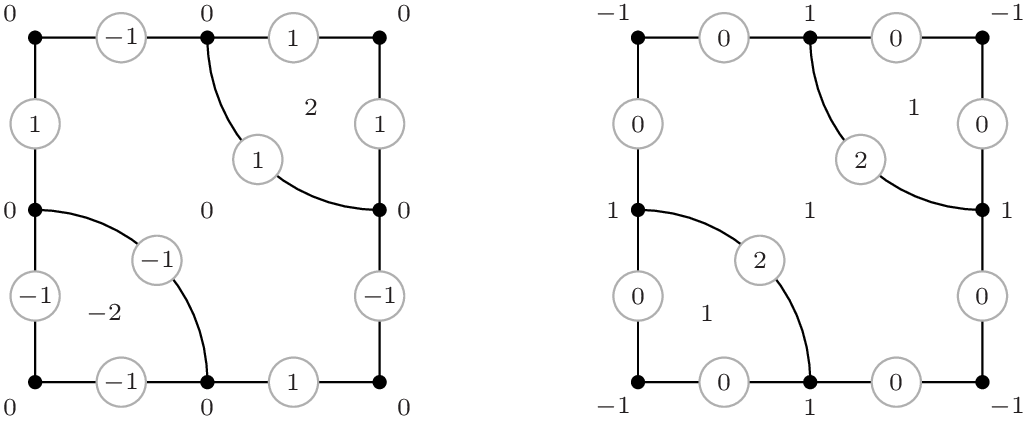}
\caption{The correction terms
$\Delta\sigma$ (left) and $\Delta\eta$ (right)
in \autoref{th.main}. See \autoref{rem.maim} for further details on this figure.}
\label{fig.main}
\endfigure
\endtheorem

It is worth emphasizing that, in both statements of \autoref{th.main}, the
knots $K^*$ contribute \emph{only} through the slopes: all other terms
depend on the links $L^*$ only.

\begin{remark}\label{rem.maim}
\autoref{fig.main} shows the correction terms
$\Delta\sigma(\slope',\slope'')$ (left)
and $\Delta\eta(\slope',\slope'')$ (right) in \autoref{th.main}.
The domain is the square $[-\infty,\infty]^2\ni(\slope',\slope'')$,
and the curve in the figures is the hyperbola $\slope'\slope''=1$.
For the term $\Delta\Gs$, we have an explicit formula, which was
actually found in the course of the proof:
\[
\Delta\Gs(\slope',\slope'')=\sg\slope'-\sg\Bigl(\frac{1}{\slope'}-\slope''\Bigr)
\label{eq.Delta.sigma}
\]
(see \autoref{cor.as.slope} for the conventions on~$\sg$);
in spite of its appearance, it is symmetric in $(\slope',\slope'')$.
Note that always
\[*
\Delta\eta(\slope',\slope'')=\ls|\Delta\Gs(\slope',\slope'')|\pm1.
\]
Intuitively, this means that the matrices always differ by an extra
eigenvector rather than by the eigenvalue of a common eigenvector.
\end{remark}

\subsection{Proof of \autoref{th.main}} \label{section.proof}
In
view of our uniform conventions on the signature, nullity, and slope (see
\autoref{rem.extension} and the respective definitions),
we can start with patching the components of~$L^*$ on which
$\Go^*$ vanish. Thus, from now on, we assume both characters nonvanishing.

If $L^*$ becomes empty ($\Go^*=1$) or was empty in the first place, we
take for~$L^*$ a small unknot contained in a ball disjoint from~$K^*$,
endowed
with its own color and arbitrary nonvanishing unitary character; this
change does not affect any of the quantities involved. (Alternatively, one
can also repeat the computation below,
taking into account the difference between the
slopes of $H_{1,0}$ and $H_{1,m}$, $m>0$, see \autoref{ex.Hopf}
and \cf. the proof of \autoref{lem.spe}.)

Thus, assume that $L^*\ne\varnothing$ and let
$(N^*,D^* \cup F^*)$ be a \emph{special} spanning pair for
$(\mathbb{S}^*,K^* \cup L^*)$, \ie, such that
$D$ is a disk and $F\cap D\ne\varnothing$. (The
existence of special spanning pairs follows from \cite[Lemma 4.1]{DFL} and
the obvious fact that, if $F\ne\varnothing$,
one can always create an extra pair of intersection
points.)
We will construct an appropriate spanning pair
 for $(\mathbb S,L)$ by cut and paste on the manifolds $(N^*,D^* \cup F^*)$.

\begin{figure}\centering
\resizebox{\linewidth}{!}{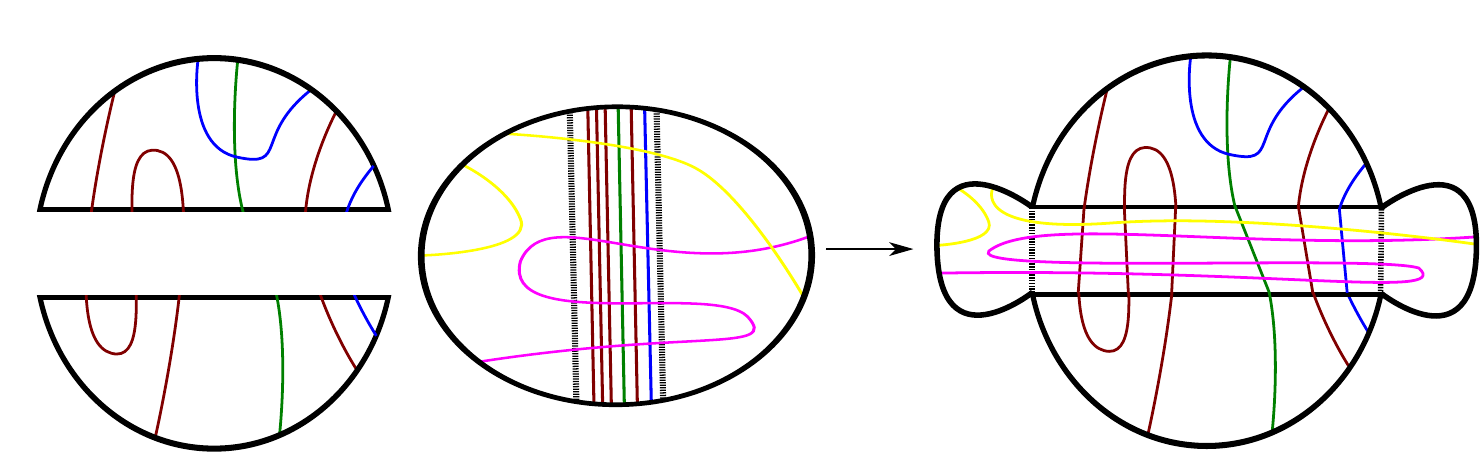}
\caption{The third diagram represents the pair $(N,F)$ used to
compute the signature of the splice of $K'\cup L'$ and $K''\cup L''$.
This pair is obtained identifying parts of the boundary of
$(N'\smallsetminus B', F' \cap ( N'\smallsetminus B'))$ and $(N'',\mathcal{D}''\cup F''$).}
\label{fig.DFL}
\end{figure}
The main idea is as in~\cite{DFL},
see \cite[Figure 4]{DFL} (reproduced here as \autoref{fig.DFL} for the
reader's convenience):
we cut off tubular neighborhoods of~$D^*$
(thus, passing to the literal extension
$\tilde\Gs_{K^*\cup L^*}(1,\Go^*)=\Gs_{L^*}(\Go^*)$, see \autoref{lem.spe}),
fill one
of the gaps with a ``standard'' spanning pair for the generalized Hopf link,
and attach the second manifold to the result. (If it
is \autoref{th.DFL} that is to be proved, the two tubular neighborhoods have
already been removed in the definition and the first step is skipped.)
An
important difference is the fact that, while in~\cite{DFL} we kept the
ramification surfaces (the colored curves in \autoref{fig.DFL})
inside the $4$-manifolds, here we need to carve them
out and work with $N^*\sminus(D^*\cup F^*)$, which makes the description of
the boundary more involved.
However, the extra boundary
parts acquired are $(\Go',\Go'')$-invisible $3$-manifolds
and this change does not
affect the computation of the homology groups, see \autoref{cor.plumbing}.

Thus,
fix neighborhoods $T_{F^*}$ and $T_{D^*}=B^2\times D^*$ and let

$W_{D^*\cup F^*} = N^* \sminus(T_{F^*} \cup T_{D^*})$.
The intersection
$B_1^* \cup \cdots\cup B_{m^*}^*= F^* \cap T_{D^*}$
is a
collection of parallel disks in
 $T_{D^*}$, and their boundary is the unlink $U^*\subset S^3=\partial T_{D^*}$
 with $m^*>0$
 components.
Consider the generalized Hopf link $H_{m',m''}= U' \cup U''$
and orient its components according to
the signs of the intersection points $F^* \cap D^*$.
Let $B=B^2 \times B^2$ be a $4$-ball,
and let $C\subset B$ be two families of parallel disks,
transversal to each other, so that $(B,C)$ is a spanning pair for
$(S^3, H_{m',m''})$.
Then,
\[*
(N,F):=\bigl(W_{D'} \cup B \cup W_{D''},
(F' \sminus (F'\cap T_{D'})) \cup C \cup  (F'' \sminus (F''\cap T_{D''}))\bigr),
\]
after smoothing the corners, is a spanning pair for $(\mathbb{S},L)$
(see \autoref{fig.DFL}).

The components of $F^*$ induce a
$(\mu'+\mu'')$-coloring on $H_{m',m''}$, and
the characters $\omega^*$ induce a unique character
$\Go:=(\omega',\omega'')$ on $H_{m',m''}$.
In the rest of the proof, we switch to~$\Go$, as all other characters are
essentially its restrictions.
To compute $\sigma_L(\Go)$, we consider the manifold
$$
W:= W_{D' \cup F'} \cup (B \sminus T_C) \cup W_{D'' \cup F''},
$$
and apply Wall's theorem to the first two terms in the above decomposition,
that is, to the manifold
$$
W_{1}:= W_{D' \cup F'} \cup (B \sminus T_C),
$$
for which the gluing takes place along $S^{3}\sminus H_{1,m'}$.

In the notation of \autoref{th.wall}, we have
$T= \partial T_{K'} \sqcup \, \partial T_{U'}$
and
\[*
X_{0}=S^3 \sminus H_{1,m'},\quad
X_{1}=(\mathbb{S}' \sminus  T_{K' \cup L'})\cup\invsbl{\partial T_{F'}\cap(N'\sminus T_{D'})},\quad
X_{2}=(S^{3}\sminus H_{m'',1})\cup\invsbl{\partial T_{C}},
\]%
where, as usual, embraced are $\Go$-invisible
parts.
Since $\omega'$ is admissible and nonvanishing, we also have
$H_1(T;\CC(\Go))=H_1(\partial T_{K'};\CC(1))=\CC^2$,
and, using \autoref{cor.plumbing} to ignore the $\Go$-invisible parts, we
obtain the following expressions for the spaces
$A_i\subset H_1(X_i;\C(\Go))$ in
\autoref{th.wall}:
\begin{alignat*}2
H_1(X_0;\CC(\omega))&= H_1(S^3 \sminus H_{1,m'}; \CC(1,\omega')),\qquad
 & A_0&= \langle m_{K'} \rangle,\\
H_1(X_1;\CC(\omega))&= H_1(\mathbb{S}' \sminus  T_{K' \cup L'}; \CC(\omega)),
 & A_1&= \langle a' m_{K'} + b' l_{K'} \rangle, \\
H_1(X_2;\CC(\omega))&= H_1(S^{3}\sminus H_{m'',1};\CC(\Go'',1)),
 & A_2&=\langle m_{K''} \rangle= \langle l_{K'} \rangle,
\end{alignat*}
where $-a'/b'=\slope'$ and $A_0$ and $A_{2}$ are computed using \autoref{ex.Hopf}.
In accordance with the settings of \autoref{s.additivity}, the orientation on
$\partial T_{K'}$ is induced from $\sphere'\sminus K'$;
hence, $m_{K'}\circ l_{K'}=-1$ in~$T$
and \autoref{cor.as.slope} yields
\[*
\sign f=\sgn(\infty,\slope',0)=-\sg\slope'.
\]
Applying \autoref{th.wall} and subtracting $\sign(N')$,
we arrive at
\[*
\begin{alignedat}3
\sign^\Go(W_{1})-\sign(N')
 &=\bigl[\sign^\Go(W_{D' \cup F'})-\sign(N')\bigr]
  &&+ \sign^\Go (B \sminus T_C)
  &&+\sg\slope'\\
 &=\sigma_{L'}(\omega')
  &&+ \delta_{\lambda'}(\Go') \delta_{\lambda''}(\Go'')
  &&+ \sg\slope'
\end{alignedat}
\]
(see \autoref{lem.spe} for the first term,
which is $\tilde\Gs_{K'\cup L'}(1,\Go')$, and
\cite[Lemma 4.2]{DFL} for the second term,
which, up to the summand $-\sign(B)=0$, is the signature of $H_{m',m''}$.)

Now, we use \autoref{th.wall} again, this time for the decomposition
 $W = W_{1} \cup W_{D'' \cup F''}$.
Since $\omega''$ is also admissible and nonvanishing,
we have
$H_1(T; \CC(\Go))  = H_1(\partial T_{K''}; \CC(1))=\CC^2$
for the new space~$T$
and,
arguing as above and ignoring the $\Go$-invisible parts, we obtain
\begin{alignat*}2
H_1(X_0;\CC(\Go))&= H_1(S^3 \sminus H_{1,m''}; \CC(1,\omega'')),\qquad
 & A_0&= \langle m_{K''} \rangle= \langle l_{K'} \rangle \\
H_1(X_1;\CC(\Go))&= H_1(\sphere' \sminus  T_{K' \cup L'}; \CC(\Go)),
 & A_1&= \langle a' m_{K'} + b' l_{K'} \rangle
       = \langle a' l_{K''} + b' m_{K''} \rangle,\\
H_1(X_2;\CC(\Go))&= H_1(\sphere''\sminus T_{K''\cup L''}; \CC(\Go)),
 & A_2&= \langle a'' m_{K''} + b'' l_{K''} \rangle,
\end{alignat*}
where $-a^*\!/b^*=\slope^*$.
Again, $m_{K'}\circ l_{K'}=-1$ and, hence, $m_{K''}\circ l_{K''}=1$,
contrary to the usual convention of \autoref{s.additivity}.
Hence, by \autoref{cor.as.slope},
\[*
\sign f=-\sgn\Bigl(\infty,\frac1{\slope'},\slope''\Bigr)
 =\sg\Bigl(\frac1{\slope'}-\slope''\Bigr),
\]
and
$$
\sign^{\Go}(W) = \sign^{\Go}(W_{1}) + \sign^{\Go}(W_{D'' \cup F''})
 -\sg\Bigl(\frac1{\slope'}-\slope''\Bigr).
$$
To complete the proof of the signature formula,
with the correction term $\Delta\Gs$ given by~\eqref{eq.Delta.sigma},
there remains to subtract the sum
$\sign(N')+\sign(N'')$
and observe that
\[*
\sign(W_{D'} \cup B \cup W_{D''})=\sign(N')+\sign(N'').
\]
For the latter statement,
one can either refer to~\cite{DFL} or
directly repeat the computation above for the ordinary signature,
when all ``slopes'' vanish.

For the nullity formula, consider the Mayer--Vietoris exact sequence related
to
\[*
X:=\Bbb S \sminus T_L=X'\cup X'',
\]
where the manifolds $X^*:= \Bbb S^* \sminus T_{K^* \cup L^*}$ are
identified along the common boundary component
$\partial_{K'} X' \cong \partial_{K''} X''$.
We have
\begin{multline*}
0 \longto
  \mathcal{K} \longto
  H_1(\partial_{K^*} X^*; \CC(1)) \longto\\
  H_1(X';\CC(\Go)) \oplus H_1(X''; \CC(\Go)) \longto
  H_1(X; \CC(\Go)) \longto
  \CC \longto 0,
\end{multline*}
where
$\mathcal{K}=\mathcal{Z}_{K' \cup L'}(\omega') \cap \mathcal{Z}_{K'' \cup L''}(\omega'')$.
Hence,
$$
\eta_L(\omega',\omega'')
 =\tilde\eta_{K' \cup L'}(1, \omega') + \tilde\eta_{K'' \cup L''}(1,\omega'')
 +\dim\mathcal{K}-1.
$$
Clearly,
$\dim\mathcal{K}=1$ if
$\KK_{K'\cup L'}(\Go')=\KK_{K''\cup L''}(\Go'')$, and $\dim\mathcal{K}=0$
otherwise.
Since $m_{K'}=l_{K''}$ and $l_{K'}=m_{K''}$ in the homology of
$\partial_{K'} X' \cong \partial_{K''} X''$,
we have $\dim\mathcal{K}=1$ if and only if $\slope'=1/\slope''$,
and there remains to apply \autoref{lem.spe} to relate
$\tilde\eta_{K^*\cup L^*}(1,\omega^*)$ and $\eta_{L^*}(\Go^*)$.
A case-by-case analysis gives us \autoref{fig.main}, right, for which we
could not find a ``nice'' formula.
\qed

\section{Skein relations for the signature}\label{S.skein}

We
conclude the paper with another illustration, \viz. we develop the concept
of slope for tangles with four
marked loose ends and analyze its relation to the signature.

For any tangle $T$ with four fixed ends in a homology $3$-ball $\mathbb{B}$
and any \emph{generic} character $\omega$ in $(\CC^\times \sminus 1)^\mu$, we
 define the slope $\kappa_T(\omega)\in\CC \cup \infty$. We show that it
 can be computed as
the quotient of the Conway polynomials of the links obtained by
 attaching to~$T$ certain \emph{elementary tangles}.
Then, we
define the sum of tangles and
prove that the signatures of
 three pairwise sums of three tangles are related by
the sign (as in \autoref{cor.as.slope}) of
their three slopes.
Finally, we use these results to derive the conventional skein relations
for the signature as in~\cite{CF}.

\subsection{Preliminaries}

Let $\ball$ be a homology $3$-ball with boundary
$S:=\partial\ball\cong S^2$.
We fix an oriented equator $E\subset S$ and four points $A_1,\ldots,A_4\in E$
ordered according to the orientation of~$E$. (We number these points
cyclically, so that $A_{i+4}=A_i$. In other words, the index takes values in
$\Z/4$.)

\definition
A \emph{tangle} is a smooth
compact oriented submanifold $T\subset\ball$ of dimension~$1$ such
that
\roster*
\item
$T$ is transversal to the boundary~$S$;
\item
$T\cap S=\{A_1,A_2,A_3,A_4\}$, and
\item
$T$ has incoming branches at $A_1,A_2$ and outgoing branches at $A_3,A_4$.
\endroster
Similarly to links, a $\mu$-\emph{coloring} on $T$ is a surjective function
$\pi_0(T) \rightarrow \{1,\dots,\mu\}$.
\enddefinition

Two tangles $T$ and $T'$ are equivalent if there exists
an orientation preserving homeomorphism
$(\ball,T)\to(\ball',T')$ taking $E$ to~$E'$ (respecting the orientation)
and $A_i$ to~$A_i'$, $i\in\Z/4$.

Given a tangle $T\subset\ball$, we let
$\ball_T=\ball\sminus \Tub T$ and
$S^\circ:=S\sminus \{ A_1,A_2,A_3,A_4\}$.
Let~$n$ be the number of components of~$T$. Then, we have
$H_1(\ball_T)\cong\Z^n$, and a
character $\Go\:H_1(\ball_T)\to\C\units$ is determined by its
values~$\Go_i \neq 1$ on the meridians~$m_i$ about the components
$T_i\subset T$. If $T$ is colored, we assume that $\Go$ takes equal values on
the meridians of the components having the same color.
We denote by $\partial\Go$ the restriction of~$\Go$ to the boundary
sphere~$S^\circ$. It is uniquely determined by the restrictions
$\Go[i]$
of $\Go$ to the meridian about $A_i$, $i \in \ZZ/4$.
Obviously,
$\Go[1]\Go[3]=\Go[2]\Go[4]$. In fact, each of $\Go[1]$, $\Go[2]$ equals one of
$\Go[3]$, $\Go[4]$ and \latin{vice versa}.

\example\label{ex.tau+-}
\begin{figure}
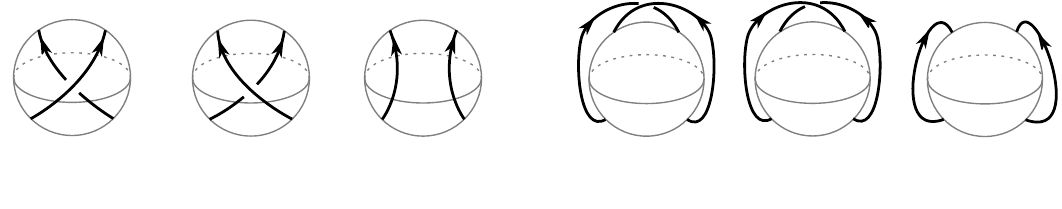
\caption{Two views of the basic tangles $\tau_{+},\tau_{-}$ and $\tau_{0}$.}
\label{f.T+-}
\end{figure}
The \emph{basic tangles}
$\tau_\pm,\tau_0\subset B^3$ are the three tangles shown in \autoref{f.T+-}
(``inside'' and ``outside'' views).
These tangles can be described as the intersection of a link undergoing a skein
transformation
with a small neighborhood of
the relevant crossing of the diagram.
\endexample

\convention\label{conv.character}
We assume that two branches of the tangle connect $A_1$ to $A_3$ and $A_2$ to $A_4$.
If the tangle is colored, we denote the two colors (that may
coincide) assigned to theses branches by~$-$ and~$+$. Therefore, we assume
that
\[*
 \Go_-:=\Go[2]=\Go[4]\ne1, \qquad \Go_+:=\Go[1]=\Go[3]\ne1.
 \]
The other convention  $\Go[1]=\Go[4]\ne1$ and $\Go[2]=\Go[3]\ne1$
can be treated similarly.
\endconvention

\begin{figure}
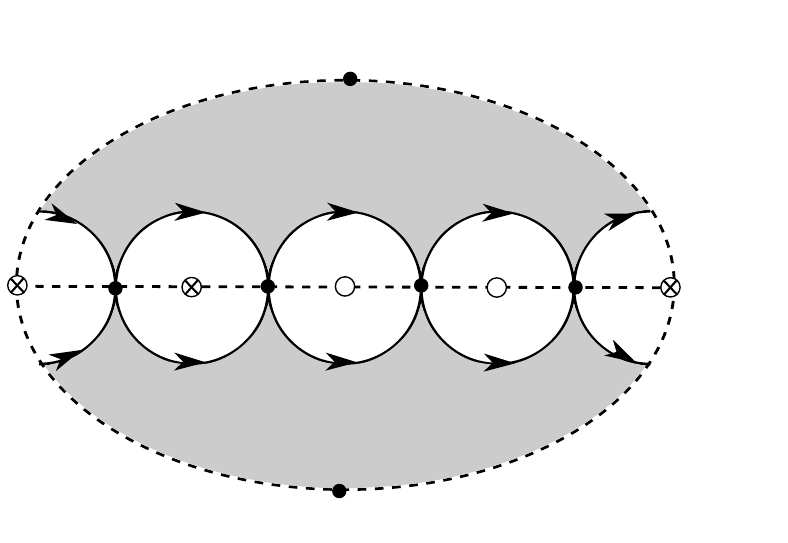
\caption{The \CW-decomposition of $S^\circ$.}
\label{fig.sphere}
\end{figure}

Up to homotopy equivalence, the punctured sphere $S^\circ$ has a
\CW-decomposition shown in \autoref{fig.sphere}
(where the sphere is cut along a ``geographic'' meridian):
there are two $2$-cells
(about each of the poles), eight $1$-cells $a_i$, $b_i$, and four
$0$-cells $B_i$, $i\in\Z/4$.
 We choose the lifts of the cells so that, for
$i\in\Z/4$, with coefficients in $\C(\Go)$,
\[*
\partial a_i=B_i-B_{i-1},\qquad
\partial b_i=\Go[i]^{\epsilon_i}B_i-B_{i-1},
\]
where $\epsilon_1=\epsilon_2=-1$ and $\epsilon_3=\epsilon_4=1$.
We define $$\tilde{c}_- :=  b_1 +  \Go_+^{-1} a_2 + \Go_+^{-1} b_3 + a_4 \qquad \tilde{c}_+:=  a_1 +  b_2 + \Go_-^{-1} a_3 + \Go_-^{-1} b_4. $$

\lemma \label{lem.basis}
For any
unitary character $\omega$ as in \autoref{conv.character}, the cycles
$$
c_-:=(1-\Go^{-1}_-)^{-1} \cdot \tilde{c}_- \qquad
c_+:=(1-\Go^{-1}_+)^{-1} \cdot \tilde{c}_+
$$
form a standard symplectic basis of  $H_1(S^\circ;\C(\Go))=\CC^2$,
in the sense of \autoref{s.additivity}.
Furthermore, the chain $c_\pm$ bounds a $2$-cell in $\ball_{\tau_\pm}$
\rom(see \autoref{ex.tau+-} and \autoref{f.T+-}\rom).
\endlemma

Note that $(c_-,c_+)$ is a basis for $H_1(S^\circ;\C(\Go))$, and
will be used
as such, for \emph{any}
character~$\Go$ as in \autoref{conv.character}. However, it is a
\emph{symplectic} basis only if $\Go$ is unitary, as otherwise we do not have
a well-defined intersection index.

\proof[Proof of \autoref{lem.basis}]
We compute the intersections in the maximal abelian covering of~$S^\circ$
and specialize the result at~$\Go$.
Since $\Go_\pm\ne1$ and the space~$S^\circ$
has homotopy type of a wedge of three circles,
we have
\[*
\dim H_1(S^\circ;\C(\Go))=2.
\]
The oriented loops $\tilde{c}_-$ and $\tilde{c}_+$ could also be defined as
the lifts of the cells of $S^\circ$:
$$
b_1+a_2+b_3+a_4 \qquad a_1+b_2+a_3+b_4
$$
starting at the chosen lift of $B_4$.
They are in general position and intersect in four points (with signs according to the orientations):
$$
(b_1+a_2+b_3+a_4) \cap (a_1+b_2+a_3+b_4) = \{ -B_4, B_1, -B_2, B_3 \}.
$$
Since  $\omega (a_1 - b_1)= \omega_+^{-1}$, $\omega (a_2 - b_2)= \omega_-^{-1}$,
and $\omega (a_3 - b_3)= \omega_+$,
we obtain
$$
\tilde c_-\circ\tilde c_+
 = -1 + \omega_+^{-1} - \omega_+^{-1} \omega_-^{-1} + \omega_-^{-1}=-(1-\omega_-^{-1})(1-\omega_+^{-1}).
$$
Hence
$c_-\circ c_+=-1$.
Since $\tilde c_-$ and~$\tilde c_+$ are represented
in the covering by honest loops, which
can be moved off, both classes and their multiples $c_-$, $c_+$ are
isotropic; thus, the latter
constitute a standard symplectic basis of
$H_1(S^\circ;\C(\Go))=\CC^2$.

 Using the relations $a_1 + a_2 + a_3 + a_4=b_1 +b_2 + b_3 +b_4=0$,
 one sees that $\tilde{c}_\pm$ can be defined as the lift of
 $h_\pm (a_4 + b_4 + a_1 + b_1)$ where $h_\pm$ is a positive or negative
 (from the south pole) half twist exchanging $A_1$ and $A_2$. This shows that $\tilde{c}_\pm$ bounds a $2$-cell in $\ball_{\tau_\pm}$.
\endproof

\subsection{Slope of a tangle}
Fix a tangle $T\subset\ball$ and consider a
character~$\Go$
as in \autoref{conv.character}.
Then, by \autoref{cor.plumbing} and \autoref{lem.basis},
there are canonical isomorphisms
$$
H_1(\partial \ball_T;\C(\Go)) = H_1(S^\circ;\CC(\Go))=\C^2
$$
and,
 as in \autoref{s.defs}, we can introduce the subspace
\[*
\KK(\Go):=\KK_T(\Go):=
\Ker\bigl[\operatorname{inclusion}_*\:H_1(S^\circ;\C(\Go))\to H_1(\ball_T;\C(\Go))\bigr].
\]

\definition
Suppose that
$\Go$ is nonvanishing and
$\dim\KK(\omega)=1$, \emph{i.e.} $\KK(\Go)$ is generated by a vector $u_-c_-+u_+c_+ $. Then the
\emph{slope} of $T$ at $\omega$ is the quotient
\[*
\kappa_T(\Go):=-\frac{u_-}{u_+}\in\C\cup\infty.
\]
As before, we extend~$\kappa_T$ to all characters by patching the
components of the tangle on whose meridians $\Go$ takes value~$1$.
Recall that we \emph{always} assume
that $\Go_\pm\ne1$ (see \autoref{conv.character}).
\enddefinition

Propositions~\ref{prop.symmetry} and~\ref{prop.unit}
and their proofs extend to $\kappa_T$
literally; in particular, if $\omega$ is a unitary character, then
$\kappa_T(\Go)$ is well-defined and real (possibly infinite).

The \emph{sum} $T'\tsum T''$
of two tangles $T'\subset\ball'$ and
$T''\subset\ball''$ is the link $T'\cup{-T''}$ in the homology sphere
$\sphere:=\ball'\cup_\partial -\ball''$, where the attaching homeomorphism
$\Gf\:S'\to S''$ restricts to an orientation preserving homeomorphism
$E'\to E''$ taking $A'_i$ to $A''_i$ for each $i=1,\ldots,4$.
If the tangles are colored,
we also ask the gluing to respect the colors~$-$ and~$+$.
Then, the result is a colored link in $\mathbb{S}$,
and $-$ and $+$ become ordinary colors.

\example\label{ex.L+-}
Consider a link in the sphere~$S^3$ and denote by~$L$ the tangle obtained
by removing a small neighborhood of a crossing of its diagrams.
Then, the sums $L_*:=L\tsum \tau_*$ with the basic tangles
(see \autoref{ex.tau+-}) are the usual links $L_\pm$ and $L_0$ involved in
a skein relation.
\endexample

\theorem\label{th.tangle}
Fix a tangle $T\subset\ball$ and consider a character~$\Go$
as in \autoref{conv.character}, so that $\Go$ extends to the links
$L_\pm:=T\tsum \tau_\pm$ \rom(see \autoref{ex.L+-}\rom). If the polynomials
$\nabla_{L_\pm}(\sqrt\Go)$ do not vanish simultaneously, then the slope
$\kappa_T(\Go)$ is well defined and one has
\[*
\kappa_T(\Go)= \frac{\nabla_{0}(\sqrt{\rGo+})}{\nabla_{0}(\sqrt{\rGo-})}
\cdot \frac{\nabla_{L_+}(\sqrt\Go)}{\nabla_{L_-}(\sqrt\Go)}\in\C\cup\infty,
\]
where $\nabla_{0}(t)=(t-t^{-1})^{-1}$ is the Conway potential of the unknot.
\endtheorem
Similarly to \autoref{th.Torres}, the statement is inconclusive if
$\nabla_{L_-}(\sqrt\Go)=\nabla_{L_+}(\sqrt\Go)=0$.

\example\label{ex.tau.slope}
The
sums of the basic tangles
(see \autoref{ex.tau+-} and \autoref{f.T+-}) are as follows:
\roster*
\item
$\tau_- \tsum \tau_-$ and $\tau_+ \tsum \tau_+$ are the Hopf links (up to orientation),
\item
$\tau_- \tsum \tau_+=\tau_+ \tsum \tau_-$ is the trivial link with two components,
\item
$\tau_- \tsum \tau_0 = \tau_+ \tsum \tau_0$ is the unknot.
\endroster
Hence, by \autoref{th.tangle}, at any \emph{nonvanishing} character, we have
$$ \kappa_{\tau_-}= \infty, \qquad \kappa_{\tau_+}=0, \qquad \kappa_{\tau_0}=1,$$
where the latter slope is defined only on the diagonal $\omega_-=\omega_+$.
\endexample

\proof[Proof of \autoref{th.tangle}]
As in \autoref{proof.Torres}, assume that the character~$\Go$ nonvanishing.
A crucial observation is the fact that, up to homotopy
equivalence, the complement of
$L_\pm$ is obtained from $\ball_T$ by attaching a
single $2$-cell~$k_\pm$ along the cycle~$\tilde{c}_\pm$ (see~\autoref{lem.basis}).
We proceed exactly as in the proof of \autoref{th.Torres},
computing and comparing the torsions
of the two links~$L_\pm$.
Let $X_\pm$ be the complement of $\Tub L_\pm$ in $\mathbb{S}=\ball \cup B^3$.
Up to homotopy equivalence, the $CW$-decomposition of
$X_\pm$ is partially given as follows
$$
C_0=\{B_1,B_2,B_3,B_4\},\quad
\{\tilde c_-,\tilde c_+,a_1,b_1,a_2,b_2,a_3,b_3,a_4,b_4,\dots\}\subset C_1,\quad
\{k_{\pm},\dots \} \subset C_2.
$$
Then, we fix the following bases for $C_*(X_{\pm};\CC(\omega))$:
$$
c_0:=B_1,B_2,B_3,B_4,\quad
c_1:=\tilde c_-, \tilde c_+,\ldots,a_1-b_1,a_2-b_2,a_3-b_3,a_4-b_4,\quad
c_2:=k_\pm,\ldots.
$$
(These bases are not quite as in the definition of the torsion. However,
they differ from cellular ones by a transition matrix in~$C_1$ that is common
to~$X_-$ and~$X_+$. Thus, our choice would not affect the \emph{ratio} of the
two torsions.)

The vectors
$$
\partial_0(a_i-b_i)=(1-\omega[i]^{\epsilon_i}) B_i, \quad i=1,\dots,4,
$$
constitute a
basis for~$C_0$. The complex $C_*(X_{\pm};\CC(\omega))$ is acyclic if and only
if
\[*
\partial_1(c_2),a_1-b_1,a_2-b_2,a_3-b_3,a_4-b_4
\]
is a basis for~$C_1$.
The homomorphisms~$\partial_1$ in~$X_-$ and~$X_+$ are given by certain
matrices of the form
$$
\bmatrix
1 & 0 &\bold0\\
\ba&\bb&M
\endbmatrix,\quad
\bmatrix
0 & 1 &\bold0\\
\ba&\bb&M
\endbmatrix,
$$
respectively,
where $M$ is a certain matrix and $\ba$, $\bb$,  are
column vectors, also common for~$X_-$ and~$X_+$.
Thus, the two complexes are acyclic if and only if neither of the two
determinants vanish, and then we have
$$
\frac{\tau_{\varphi}(X_+)(\omega)}{\tau_{\varphi}(X_-)(\omega)}=-\frac{\det[\ba\,|\,M]}{ \det[\bb\,|\,M]}.
$$
Now, a linear combination
$(1-\omega_-^{-1})^{-1} u_- \tilde c_-+ (1-\omega_+^{-1})^{-1} u_+ \tilde c_+ $
represents a class in $\KK(\Go)$ if and
only if it is in the image of~$\partial_1$ in the complex
$C_*(\ball_T;\CC(\Go))$, which is given by $[\ba\,|\,\bb\,|\,M]$.
If the slope~$\kappa$ is finite,
we are looking for a vector of the form
$-\tilde\kappa \tilde c_-+\tilde c_+\in\Im\partial_1$, where
\[*
\tilde\kappa:=\frac{(1-\omega_+^{-1})}{(1-\omega_-^{-1})} \kappa
 =\frac{\nabla_{0}(\sqrt{\rGo-})}{\nabla_{0}(\sqrt{\rGo+})}\kappa.
\]
Clearly, if such a vector exists, then
$\tilde\kappa=-\bx\cdot\ba$, where $\bx$ is a solution to the
linear system
$$
\bx\cdot[\bb\,|\,M]=[1\,|\,\bold0].
$$
If $\tau_\Gf(X_-)(\Go)\ne0$, then $\det[\bb\,|\,M]\ne0$ and
$\tilde\kappa$ can be computed by Cramer's rule:
$$
\tilde\kappa
 =-[1\,|\,\bold0]\cdot[\bb\,|\,M]\1\cdot{\ba}
 =-\frac{\det[\ba\,|\,M]}{\det[\bb\,|\,M]}
 =\frac{\tau_{\varphi}(X_+)(\omega)}{\tau_{\varphi}(X_-)(\omega)},
$$
so that there is at most one solution. As in \autoref{proof.Torres},
replacing~$\Go$ with~$\Go\dm$, we conclude that the slope is well defined and
given by the expression above.
If $\tau_\Gf(X_-)(\Go)=0$, but $\tau_\Gf(X_+)(\Go)\ne0$,
the same argument with $X_\pm$ interchanged shows that the slope
is well defined and equal to~$\infty$.

The passage from~$\tau_\Gf$ to~$\nabla$ and verification of the sign is
immediate (\cf. \autoref{proof.Torres}). \endproof

\subsection{The skein relations}

Let $T'$ and $T''$ be a pair of colored tangles such that the sum
$T' \tsum T''$ is well-defined.
Clearly, any pair of characters $\Go'$, $\Go''$ on $T'$, $T''$ such that
$\partial\Go'=\partial\Go''$
gives rise to a character on $T'\tsum T''$;
we denote the latter by $\Go'\tsum\Go''$.

\theorem\label{th.skein}
Consider
three tangles $T^i\subset\ball^i$, $i\in\Z/3$, and three
characters $\Go^i$ on $\ball^i\sminus T^i$ as in \autoref{conv.character}
such that $\partial\Go^i=\const(i)$
(\ie, all three characters have the same restriction to the
common boundary sphere~$S$\rom).
Then, denoting $\kappa^i:=\kappa_{T^i}(\Go^i)$, one has
\[*
\sum_{i\in\Z/3}\Gs_{T^{i+1}\tsum T^i}(\Go^{i+1}\tsum\Go^i)
 =\sgn(\kappa^0,\kappa^1,\kappa^2)
\]
\rom(see \autoref{cor.as.slope} for the definition of~$\sgn$\rom).
\endtheorem

\autoref{th.skein} could be
derived directly from \cite[Theorem 1.1]{CC1}.
The short proof given below, although using essentially the same
argument,
fits better the framework of our paper.

\remark\label{rem.skein}
Changing the indices $(0,1,2)$ to
$(\text{nothing},\prime,\prime\prime)$, one can rewrite the conclusion of
\autoref{th.skein} in the following, less symmetric, form:
\[*
\Gs_{T\tsum T'}(\Go\tsum\Go')+\Gs_{T'\tsum T''}(\Go'\tsum\Go'')
 =\Gs_{T\tsum T''}(\Go\tsum\Go'')-\sgn(\kappa,\kappa',\kappa'').
\]
It is, essentially, this identity that is actually proved below.
\endremark

\proof[Proof of \autoref{th.skein}]
Since the orientation of $T''$ is reversed in the definition,
the operation $T' \tsum T''$ is skew-symmetric, \ie,
for any colored tangles $T'$, $T''$
and character $\Go'$, $\Go''$ such that the operation is well-defined,
we have
$$
\Gs_{T'\tsum T''}(\Go' \tsum \Go'')=-\Gs_{T''\tsum T'}(\Go'' \tsum \Go').
$$
This observation justifies \autoref{rem.skein} and,
in the latter form, the identity
is an immediate consequence of Wall's \autoref{th.wall}.
Indeed, if $(N_1,F_1)$ is a spanning pair for $T'\tsum T$ and $(N_2,F_2)$ is
a spanning pair for $T\tsum T''$, then a spanning pair for $T'\tsum T''$ is
$(N,F):=(N_1,F_1)\cup_\ball(N_2,F_2)$. The three characters $\Go,\Go',\Go''$ define a
common character~$\tilde\Go$ on $\ball\cup\ball'\cup\ball''$, which extends
uniquely to $N_1\sminus F_1$ and $N_2\sminus F_2$, and there remains to apply
\autoref{th.wall}. The few technical details related to removing
$\Tub F_i$, $i=1,2$,
are filled in as in
\autoref{section.proof}, using \autoref{cor.plumbing};
we leave this exercise to the reader.
\endproof

The following corollary generalizes and refines \cite[Theorem 5.1]{CF}.

\corollary\label{cor.skein}
Let $L_\pm,L_0\subset S^3$ be colored links involved into the skein relation at a
crossing of the diagram, and denote by $L\subset B^3$ the tangle obtained by
removing a small tubular neighborhood of the crossing \rom(see
\autoref{ex.L+-}\rom). Pick a unitary character $\Go\in(\sone)^\mu$
on~$L_\pm$ and, {\em contrary to the usual convention,
fix $\sqrt\Go$ so that $\Im\sqrt{\rGo\pm}>0$.} Then,
\[*
\sigma_{L_+}(\omega) - \sigma_{L_-}(\omega)
 = \sg \kappa_L(\omega)
 = \sg \left(\frac{\nabla_{L_+}(\sqrt\Go)}{\nabla_{L_-}(\sqrt\Go)}\right).
\]
If $\omega_-=\omega_+$, then also
\[*
\sigma_{L_\pm}(\omega) - \sigma_{L_0}(\omega)
 = \sg (\kappa_L(\omega)^{\mp 1}-1)
 = \pm \sg \left(i\cdot \frac{\nabla_{L_\pm}(\sqrt\Go)}{\nabla_{L_0}(\sqrt\Go)}\right).
\]
In both case, the second expression makes sense if at least one of the two
Conway potentials does not vanish\rom; in this case, we assert, in
particular, that the argument of $\sg$ is real.
\endcorollary

\proof
Let $T=L$, $T''=\tau_+$, and $T'=\tau_-$ or $\tau_0$ in
\autoref{th.skein}.
Since
$\tau_+ \tsum \tau_-$ and $\tau_\pm\tsum\tau_0$ are trivial links/knots
(see \autoref{ex.tau+-}),
their signature is $0$ and, hence, letting $\kappa:=\kappa_T(\omega)$,
from
\autoref{th.skein}
and \autoref{ex.tau.slope} we have
\begin{alignat*}2
\sigma_{L_+}(\omega) - \sigma_{L_-}(\omega)
 &=\sgn(\kappa,\infty,0)&&=\sg\kappa,\\
\sigma_{L_+}(\omega) - \sigma_{L_0}(\omega)
 &=\sgn(\kappa,1,0)=\sg(\kappa(1-\kappa))&&=\sg(\kappa^{-1}-1).
\end{alignat*}
To relate these expressions to the Conway potentials, we use
\autoref{th.tangle} and the following simple observation: if $\xi\in S^1$,
the difference $\xi-\xi\1=2i\Im\xi$ makes a predictable
contribution to the sign. In view of our choice of
$\sqrt{\rGo\pm}$,
this
completes the proof of the first formula. For the second one, we employ the
classical skein relation (letting
$\xi:=\sqrt{\rGo+}=\sqrt{\rGo-}$)
\[*
\nabla_{L_+}(\sqrt\Go)-\nabla_{L_-}(\sqrt\Go)
 =(\xi - \xi\1)\nabla_{L_0}(\sqrt\Go),
\]
which implies
\[*
\kappa\1-1
 = \frac{\nabla_{L_-}(\sqrt\Go)}{\nabla_{L_+}(\sqrt\Go)} - 1
 = \frac{\nabla_{L_-}(\sqrt\Go) - \nabla_{L_+}(\sqrt\Go)}{\nabla_{L_+}(\sqrt\Go)}
 = (\xi\1 - \xi)\frac{\nabla_{L_0}(\sqrt\Go)}{\nabla_{L_+}(\sqrt\Go)}
 \in\RR\cup\infty
\]
and, since $\Im\xi>0$,
\[*
\sg(\kappa\1-1)
 = \sg\left( -i \frac{\nabla_{L_0}(\sqrt\Go)}{\nabla_{L_+}(\sqrt\Go)} \right)
 = \sg\left( i \frac{\nabla_{L_+}(\sqrt\Go)}{\nabla_{L_0}(\sqrt\Go)} \right).
\]
The computation for $\sigma_{L_-}(\omega) - \sigma_{L_0}(\omega)$ is similar.
\endproof

\let\.\DOTaccent
\def\cprime{$'$}
\bibliographystyle{amsalpha}
\bibliography{bibliosplice}

\end{document}